%% file: EntropyStableCFN.tex
\documentclass[final,onefignum,onetabnum,letterpaper]{siamart220329}

\input{EntropyStableCFN_shared.tex}

\ifpdf
\hypersetup{
  pdftitle={Entropy},
  pdfauthor={Tongtong Li}
}
\fi

\begin{document}

\maketitle

\begin{abstract}
	\input{0_abstract}
\end{abstract}





\input{1_introduction} 
\input{2_preliminaries}

\input{3_CFNnetwork} 
\input{4_Design}
\input{5_Numerical}
\input{6_summary}

\section*{Acknowledgements} 
This work was supported by  DoD ONR MURI grant \#N00014-20-1-2595 (all), AFOSR grant \#F9550-22-1-0411 (AG), and DOE ASCR grant \#DE-ACO5-000R22725 (AG). 

\appendix 

\small
\bibliographystyle{siamplain}
\bibliography{reference}

\end{document}

%% file: EntropyStableCFN_shared.tex
\usepackage{soul}
\usepackage{lipsum}
\usepackage{amsfonts}
\usepackage{epstopdf}
\usepackage{epsfig}
\ifpdf
  \DeclareGraphicsExtensions{.eps,.pdf,.png,.jpg}
\else
  \DeclareGraphicsExtensions{.eps}
\fi

\usepackage{datetime}
\newdateformat{monthyeardate}{%
	\monthname[\THEMONTH] \THEDAY, \THEYEAR}

\usepackage{academicons}
\usepackage{xcolor}

\usepackage{amsmath}
\usepackage{amssymb}
\usepackage{commath}
\usepackage{mathtools}
\usepackage{bbm}
\usepackage{bm}
\usepackage{upgreek}
\usepackage{float}
\allowdisplaybreaks

\usepackage{color}
\usepackage{graphicx}
\usepackage[small]{caption}
\usepackage{subcaption}

\usepackage{relsize}
\usepackage{adjustbox}
\usepackage{algorithm}
\usepackage{algpseudocode}
\usepackage{booktabs}
\usepackage{tikz}
\usepackage{pifont}
\usetikzlibrary{bayesnet} 

\usepackage{enumitem}
\setlist[enumerate]{leftmargin=.5in}
\setlist[itemize]{leftmargin=.5in}

\newsiamremark{remark}{Remark}
\crefname{hypothesis}{Hypothesis}{Hypotheses}

\headers{Entropy stable conservative flux form neural networks}{Lizuo Liu, Tongtong Li, Anne Gelb, and Yoonsang Lee}

\title{Entropy stable conservative flux form neural networks
\thanks{
\monthyeardate\today 
}}

\author{
Lizuo Liu\thanks{Department of Mathematics, Dartmouth College, Hanover, NH 03755, USA 
(\email{lizuo.liu@dartmouth.edu}, 
\email{Anne.E.Gelb@Dartmouth.edu}, 
and \email{Yoonsang.Lee@dartmouth.edu}
)}
\and 
Tongtong Li\footnotemark[2]\thanks{Department of Mathematics and Statistics, University of Maryland, Baltimore County, Baltimore, MD 21250, USA (\email{tongtong.li@umbc.edu}
) } 
\and 
Anne Gelb\footnotemark[2]
\and 
Yoonsang Lee\footnotemark[2]
}

\usepackage{amsopn}

\usepackage{comment}
\newcommand{\ds}{\displaystyle}

\usepackage[normalem]{ulem}

\newcommand{\btheta}{{\boldsymbol\theta}}

\newcommand{\bu}{\mathbf{u}}
\newcommand{\bbR}{\mathbb{R}}
\newcommand{\cF}{\mathcal{F}}
\def\wb{\overline}

\usepackage[normalem]{ulem}

\newtheorem{rem}{Remark}[section]

\numberwithin{equation}{section}
\numberwithin{figure}{section}
\numberwithin{table}{section}

\allowdisplaybreaks

\patchcmd\newpage{\vfil}{}{}{}

%% file: 0_abstract.tex
We propose an entropy-stable conservative flux form neural network (CFN) that integrates classical numerical conservation laws into a data-driven framework using the entropy-stable, second-order, and non-oscillatory Kurganov-Tadmor (KT) scheme. The proposed entropy-stable CFN uses slope limiting as a denoising mechanism, ensuring accurate predictions in both noisy and sparse observation environments, as well as in both smooth and discontinuous regions. Numerical experiments demonstrate that the entropy-stable CFN achieves both stability and conservation while maintaining accuracy over extended time domains. Furthermore, it successfully predicts shock propagation speeds in long-term simulations, {\it without} oracle knowledge of later-time profiles in the training data. 

%% file: 1_introduction.tex
\section{Introduction}
\label{sec:introduction}

Hyperbolic systems of partial differential equations (PDEs) are fundamental in modeling the dynamics of various natural and engineered systems. These systems are characterized by their ability to describe wave propagation, making them essential for understanding phenomena involving the transport of quantities such as mass, momentum, and energy. They therefore play a pivotal role in the fields of climate science, oceanography, and sea ice dynamics \cite{LeVeque02}. For example, the shallow water equations are used extensively to model ocean currents and wave dynamics, helping to predict the movement of oceanic water masses and their impact on coastal regions \cite{Shallow_water}. Euler equations are similarly important in atmospheric dynamics models. Their general utility in large-scale flow patterns and wave phenomena make them fundamental to  modeling weather and climate systems \cite{Euler}. Iceberg and sea ice dynamics  affect both oceanic and atmospheric systems, with hyperbolic PDEs playing roles in the modeling of iceberg motion \cite{iceberg} and in the elastic deformation of sea ice \cite{Girard11, Elastic-PlasticSeaIce}. It is therefore consistent to assume that {\em unknown} systems of hyperbolic PDEs govern the dynamics that describe intricate relationships between ice, ocean, and atmosphere.  This investigation demonstrates that it may be possible  to predict future behavior of these unknown systems from relevant data acquired in the past.

The rigorous mathematical framework for hyperbolic systems of PDEs has enabled the development of accurate, efficient and robust numerical solvers for these complex processes when the governing PDEs are {\em known}. It is imperative for numerical schemes to closely mirror the mathematical structures of the underlying PDEs, as approximate solutions can develop singularities such as shocks, sharp gradients, and discontinuous derivatives, even when starting from smooth initial conditions.  Due to their ability to resolve the locally complicated structures inherent in nonlinear hyperbolic PDEs, finite volume methods \cite{LeVeque02} and discontinuous Galerkin methods \cite{Hesthaven_DG} are often adopted for long term model simulations \cite{dg_wave, dg_acoustic, fvm_unstructured, fvm_weno}.   

As large datasets pertaining to important state variables (e.g.~temperature, pressure, and velocity fields) are more readily acquired, data-driven models are increasingly being used to accurately and efficiently learn the dynamics of {\em unknown} ordinary differential equation (ODE) and PDE models, as well as to predict their future behaviors. In this regard machine learning (ML) offers a compatible framework for various dynamical system prediction problems, originally in the context of ODEs \cite{dynonet, rnn_ode, cdeeponet}, with some recent extensions to PDEs \cite{neuralPDE,deeppropnet, yin2023continuous}. Analogous to the well known method-of-lines approach, most data-driven methods for learning dynamical systems construct neural networks to represent the right hand side of a semi-discrete ODE system, whose inputs are predictions at previous time steps, which is then followed by appropriate numerical time integration techniques \cite{chenPDE, Churchill_2023}. These approaches are considered as purely data-driven since no explicit spatial structural information is included. By contrast, methods that incorporate spatial derivatives such as gradients and Laplacians are designed to embed explicit spatial structural information for dynamical system discovery \cite{long2018pdenetlearningpdesdata}. 

Methods that incorporate ideas from classical numerical conservation laws into neural network frameworks have been more recently developed to predict the long term behavior of hyperbolic conservation laws \cite{charles2024learningwenoentropystable, chencfn, RoeNet}. For example, the {\em conservative flux form neural network (CFN)} introduced in \cite{chencfn} learns the dynamics of unknown hyperbolic conservation laws by leveraging a finite volume structure. Similarly, RoeNet in \cite{RoeNet} draws inspiration from the classical Roe scheme, which approximates the Riemann solver \cite{roe}. Notably, CFN and RoeNet both describe methods that extrapolate the solutions to a future time using training data obtained in some initial time interval.   

Based on the promising results achieved with CFN and RoeNet, this investigation expands on the theme of incorporating classical numerical conservation laws, specifically to include {\em entropy stability conditions}, into the neural network framework.  We note that the method in \cite{RoeNet} is reported to not ``violate entropy'', but neither formal discussion nor computational evidence is provided.  By contrast, in this investigation we intentionally consider entropy stability in the design of our neural network in the form of slope-limiting numerical methods. As a prototype we employ the second order accurate non-oscillatory Kurganov-Tadmor (KT) scheme \cite{KTscheme2000}, which uses a minmod slope limiter. We observe that much like what is provable for classical numerical methods for conservation laws \cite{Hesthaven, LeVeque92, LeVeque02},  the combination of entropy stability, higher (here second) order accuracy, and slope-limiting within the neural network design yields better prediction capabilities for data-driven methods than the original CFN in \cite{chencfn}, which takes none of these properties into account. In particular, by design the original CFN approach tends to preserve oscillations in the solution, and therefore cannot discern true  oscillatory behavior from that induced by noise, in turn potentially leading  to incorrect flux terms and non-physical solutions. By contrast the slope-limiting property of the entropy-stable CFN scheme proposed here  effectively serves as a denoising procedure since it reduces oscillations in the solution. Moreover, since its construction is inherently tied to the spatial order of accuracy, slope limiting is more robust in providing the expected order of accuracy for the solution (in smooth regions) when compared to the original CFN approach. 

Our numerical experiments demonstrate that the entropy-stable CFN yields results that are both conservative and entropy-stable. It also correctly predicts the shock propagation speed in extended time domains without introducing non-physical oscillations into the solution. Importantly, these accurate predictions can be generated even when the training period does not contain a discontinuous profile. That is, the method does not require oracle knowledge of a later time solution space. Our new method is more efficient than techniques such as RoeNet \cite{RoeNet}, since it bypasses the need for computing a matrix pseudo-inverse, which has the added risk of generating instabilities during training and testing.  Finally, our method's network architecture is simpler and more compact, further enhancing efficiency.  

The rest of this paper is organized as follows. In \Cref{sec:preliminaries} we review fundamental properties of numerical conservation laws and the entropy condition, as well as discuss three entropy-stable numerical methods that can potentially provide the building blocks for our new entropy-stable CFN, namely the KT, the Lax Wendroff and the modified Lax Wendroff schemes.  We then construct the entropy-stable CFN in \Cref{sec:entropy-stable network}. In \Cref{sec:experimentdesign} we provide details for our experimental design and the metrics used to validate our results. Computational experiments are performed in \Cref{sec: numerical}, where we also analyze our method for its sensitivity to noise, training data resolution, and choice of parameters.  Some concluding remarks and ideas for future investigations are given in \Cref{sec:conclusion}. 

%% file: 2_preliminaries.tex
\section{Conservation Laws}
\label{sec:preliminaries}

We seek to learn the dynamics of an unknown hyperbolic system of conservation laws. For self-containment purposes, we review the fundamental properties that are important in constructing numerical solvers for known governing equations.

For ease of presentation, we consider the 1D scalar conservation law of the form  
\begin{equation}
\label{eq: conservation_law}
\frac{\partial}{\partial t} u + \frac{\partial}{\partial x} f (u) = 0, \quad x \in (a,b), \quad t \in (0,T),
\end{equation}
 with appropriate boundary and initial conditions.  Here $f(u)$ is a smooth flux function of the conserved variable $u = u(x,t)$. A distinctive feature of \cref{eq: conservation_law} is the spontaneous formation of shock discontinuities, which can occur even when the initial conditions are smooth. This also makes the problem notoriously difficult to solve.  In this regard the {\em entropy condition}, which requires all associated admissible entropy pairs $(U(u), F(u))$, where $U(u)$ is the entropy function and $F(u)$ is the entropy flux, to satisfy the additional inequality
\begin{equation} 
\label{eq:entropy inequality}
\frac{\partial}{\partial t} U(u) +  \frac{\partial}{\partial x} F(u) \leq 0, 
\end{equation}
plays a critical role in developing numerical methods that ensure that the corresponding computed solution of \cref{eq: conservation_law} is physically meaningful \cite{Dafermos2000HyberbolicCL, Lax1987HyperbolicSO,  Smoller1983ShockWA}. 

\subsection{Numerical schemes for conservation laws}
\label{sec: numerical conservation laws}
Entropy-stable schemes have been developed in a variety of contexts.  This investigation employs the (modified) Lax-Wendroff and the Kurganov-Tadmor schemes, which we now describe. Both are centered schemes based on uniform spatial domain grid points and are chosen as prototypes here because of their relatively simple constructions.  Other methods may be more expedient depending on the particular problem environment.  For example, discontinuous Galerkin methods are particularly well-suited for complex geometries. Future investigations will consider other entropy-stable schemes as appropriate. 

We begin by discretizing the spatial domain with grid points $\{x_j\}_{j = 0}^n$ and the temporal domain with time instances $\{t_l\}_{l =0}^L$ as\footnote{While the modified Lax-Wendroff scheme is designed explicitly for uniform temporal and spatial discretizations, the Kurganov-Tadmor scheme uses the method of lines.  In this case we employ uniform time instances to ensure a direct comparison for our numerical results.}
\begin{equation}
\label{eq:spatialgrid}
    x_j=j \Delta x, \quad \Delta x = \frac{b-a}{n} \quad \text{and} \quad t_l = l \Delta t, \quad \Delta t = \frac{T}{L}.\end{equation} 
  We then use cell averages to approximate the solution, $\bu(t) = \{ \wb{u}_j(t) \}_{j = 0}^n$, given by
\begin{equation}
\label{eq: cell_average}
\wb{u}_j(t) = \int_{x_j - \frac{\Delta x}{2}}^{x_j + \frac{\Delta x}{2}} u(x,t)\,dx, \quad j=1, \cdots, n-1, 
\end{equation}
with appropriate boundary conditions enforced for $\wb{u}_0$ and $\wb{u}_n$.
Traditional numerical conservation laws aim to construct accurate, robust, and efficient solvers for $\wb{u}_j(t)$, $j = 0,\dots,n$, for a given flux function $f$.  To retain the proper shock speed, numerical solvers for \cref{eq: conservation_law} must be written in  conservative form \cite{Hesthaven, LeVeque92, LeVeque02}. In particular, using the corresponding semi-discretized equation of \cref{eq: conservation_law},  each cell average in \cref{eq: cell_average} is updated by calculating the flux differences at the cell edges, that is, by solving
\begin{equation}
\frac{d}{dt} \wb{u}_j + \frac{1}{\Delta x} \left(  f_{j+1/2} - f_{j-1/2} \right) = 0,
\label{eq: semi-discrete conservative scheme}
\end{equation}
where $f_{j+1/2}$ is a numerical flux approximating the true flux at the cell edge $x=x_{j+1/2}$. Entropy stability is obtained by satisfying the discrete entropy inequality corresponding to \cref{eq:entropy inequality}, and is given by 
\begin{equation}
\label{eq: entropy inequality semi}
\frac{d}{dt} U(\wb{u}_j) + \frac{1}{\Delta x} \left(  F_{j+1/2} - F_{j-1/2} \right) \leq 0.
\end{equation}
A detailed analysis of entropy stability  for difference approximations may be found in \cite{TadmorConservation2003}.

\subsubsection{Lax-Wendroff (LW) scheme} \label{sec:LxWscheme}
The Lax-Wendroff (LW) scheme \cite{L-Wscheme} is a conservative-form second-order difference scheme in both space and time, and is given for the generic scalar conservation law \cref{eq: conservation_law} as
\begin{equation}
\begin{aligned}
\wb{u}_{j}^{l+1} &= \wb{u}_{j}^{l} 
- \frac{1}{2}\frac{\Delta t}{\Delta x}
\left[ f(\wb{u}_{j+1}^{l}) - f ( \wb{u}_{j-1}^{l} )\right] 
\\ &+ \frac{1}{2}\left(\frac{\Delta t}{\Delta x}\right)^2
\left[ A_{j+1/2}\left( f(\wb{u}_{j+1}^l) - f(\wb{u}_{j}^l)\right) - A_{j-1/2}\left( f(\wb{u}_{j}^l) - f(\wb{u}_{j-1}^l)\right)  \right], \\
\end{aligned}
\label{eq: LWS}
\end{equation}
where \(A_{j\pm 1/2}\) is the Jacobian matrix evaluated at \(\frac{1}{2}\left(u_j^{l} + u_{j\pm 1}^{l}\right)\).

\subsubsection{The modified Lax-Wendroff  (modLW) scheme} \label{sec:modLWscheme}
Although widely employed for solving systems of conservation laws, it is well known that in some cases the LW scheme converges to a nonphysical weak solution \cite{Harten-nonunique-Lw, MacCormack-nonunique-lw}. The modified Lax-Wendroff (modLW) scheme  introduced in \cite{mlw} adds a nonlinear artificial viscosity term to the LW operator. In addition to maintaining conservative form and second-order accuracy, it also has the desirable property that the limiting solution of a converging solution is guaranteed to satisfy the entropy condition \cref{eq:entropy inequality}.  In 1D, the modLW scheme has the form\footnote{Although cell average notation is used in \cref{eq: cell_average}, we note that finite difference and finite volume schemes are equivalent in 1D.}
\begin{equation}
\begin{aligned}
\wb{u}_{j}^{l+1} &= \wb{u}_{j}^{l} 
- \frac{1}{2}\frac{\Delta t}{\Delta x} 
\left[ f(\wb{u}_{j+1}^{l}) - f ( \wb{u}_{j-1}^{l} )\right] 
+ \frac{1}{2} \left(\frac{\Delta t}{\Delta x}\right)^2 
\Delta_{-}\left[ \frac{\Delta_{+} f ( \wb{u}_{j}^{l})}{\Delta_{+} \wb{u}_{j}^{l}} \Delta_{+}f ( \wb{u}_{j}^{l} ) \right]\\
&+ \dfrac{\Delta t}{\Delta x} \Delta_{-} \left[ 
C \gamma\left( \frac{|\Delta_{+} \wb{u}_{j}^{l}|}{(\Delta x)^{\alpha}} \right) \Delta_{+} \frac{\partial f}{\partial u}\left( \wb{u}_{j}^{l} \right) \Delta_{+} \wb{u}_{j}^{l} \right], \\
\end{aligned}
\label{eq: modLWS}
\end{equation}
where $\wb{u}_{j}^{l}=\wb{u}_{j}(t_l)$, \(\Delta_{+} u_{j} = u_{j+1} - u_{j}, \Delta_{-} = u_{j} - u_{j-1} \), and \(\gamma(\cdot)\) is given by
   \[
   \gamma\left( s \right) = \left\{
   \begin{aligned}
       0 \quad \text{for} \quad |s| < 1,\\
       1 \quad \text{for} \quad |s| \ge 1.
    \end{aligned}
   \right.
   \]
Constant parameters $C$, $\alpha  > 0$ are chosen to guarantee conservation and entropy stability.  In our numerical experiments we use $C = 0.05$ and $\alpha = 1$. Details concerning admissible parameter choices may be found in \cite{mlw}. 

\subsubsection{Kurganov-Tadmor (KT) scheme}
\label{sec:KTscheme}
The second order entropy-stable Kurganov-Tadmor (KT) scheme in \cite{KTscheme2000} is also designed to be nonoscillatory. Specifically, it is total-variation diminishing (TVD) in 1D  and satisfies the maximum principle for  2D scalar conservation laws \cite{KTscheme2000}.
Consequently,  the introduction of nonphysical oscillations in piecewise smooth solutions may be entirely avoided. The scheme admits the conservation form
\begin{equation}
\frac{d}{d t} \wb{u}_j(t)=-\frac{H_{j+1/2}(t)-H_{j-1 / 2}(t)}{\Delta x},
\label{eq: Kurganov Tadmor Scheme}
\end{equation}
with numerical flux
\begin{equation}
H_{j+1/2}(t):=\frac{f\left(u_{j+1 / 2}^{+}(t)\right)+f\left(u_{j+1/2}^{-}(t)\right)}{2}-\frac{a_{j+1 / 2}(t)}{2}\left[u_{j+1 / 2}^{+}(t)-u_{j+1 / 2}^{-}(t)\right].
\label{eq: KT flux}
\end{equation}
Observe that \cref{eq: KT flux} requires precise information regarding local propagation speed at the cell boundary $x_{j+1/2}$, as determined by
\[
a_{j+1 / 2}(t):=\max \left\{\rho\left(\frac{\partial f}{\partial u}\left(u_{j+1 / 2}^{+}(t)\right)\right), \ \rho\left(\frac{\partial f}{\partial u}\left(u_{j+1 / 2}^{-}(t)\right)\right)\right\}. 
\]
Here the intermediate values $u_{j+1 / 2}^{ \pm}$ are defined as
\begin{equation}
u_{j+1 / 2}^{+}:=\wb{u}_{j+1}(t)-\frac{\Delta x}{2}\left(u_x\right)_{j+1}(t), \quad u_{j+1 / 2}^{-}:=\wb{u}_j(t)+\frac{\Delta x}{2}\left(u_x\right)_j(t),
\label{eq: upm}
\end{equation}
where $(u_x)_j$ is the approximate spatial derivative at grid location $x_j$. The nonoscillatory behavior of the KT scheme hinges on how this is computed:
\begin{equation}
\left( u_{x} \right)_{j} = \phi\left( r \right)\left( \wb{u}_{j + 1} - \wb{u}_{j} \right), \quad r = \dfrac{\wb{u}_{j} - \wb{u}_{j-1}}{\wb{u}_{j+1} - \wb{u}_{j}},
\end{equation}
with total-variation stability obtained via the minmod slope limiter
\begin{equation}
\label{eq: minmod limiter}
\phi(r) =  \phi_{minmod}\left( r \right) = \max\left( 0, \min\left( r, \dfrac{1+r}{2}, 1 \right) \right).
\end{equation}

\subsection{Time integration}
\label{sec:timeintegration}
As already noted, both the LW and modLW schemes are simultaneously temporally and spatially discretized, yielding a one-step time integration.  By contrast,  \cref{eq: Kurganov Tadmor Scheme} is solved using the method of lines.  To ensure a viable stability region for an efficient choice of time step $\Delta t$,  we utilize the  third-order total variation diminishing Runge-Kutta (TVDRK3) method \cite{Shu98} as the time integration technique, which can be employed on any generic time-dependent ODE given by
$$ \frac{dz}{dt}= \mathcal{L}(z),$$
where $\mathcal{L}$ is a known operator on $z$. The TVDRK3 method provided in \Cref{alg: TVDRK3} advances the solution $z^{l-1}$, $l \ge 1$, at the current time level ${t}_{l-1}$ to the solution $z^l$ at the next time level ${t}_{l}$. The value $z^0$ corresponds to the initial conditions $u(x,0)$ in \cref{eq: conservation_law}. 

\begin{algorithm}[H]
\caption{TVDRK3 time integration method for a single time step starting at time level $t^{l-1}$}\label{alg: TVDRK3}
\begin{algorithmic}
\State INPUT: ${z}^{l-1}$, $\mathcal{L}(z^{l-1})$ and $\Delta t$.
\State OUTPUT: The solution ${z}^{l}$ at time level $t_l$.
\State $ z^{(1)}=z^{l-1}+ \Delta t \, \mathcal{L}(z^{l-1})$,
\medskip
\State $ z^{(2)}=\frac{3}{4}z^{l-1} + \frac{1}{4}z^{(1)} + \frac{1}{4}\Delta t \, \mathcal{L}(z^{(1)})$,
\medskip
\State $ z^{l}=\frac{1}{3}z^{l-1} + \frac{2}{3}z^{(2)}+\frac{2}{3}\Delta t \, \mathcal{L}(z^{(2)}).$
\end{algorithmic}
\end{algorithm}

%% file: 3_CFNnetwork.tex
\section{Entropy-stable network}
\label{sec:entropy-stable network}

We now have all the ingredients needed to define our new {\em entropy-stable} network, which is described below.  Since our method extends upon the approach used in \cite{chencfn}, we will adopt the same framework as was used there.

\subsection{Problem setup}

We assume that training data for the unknown system of conservation laws are available as a set of trajectories of discrete temporal snapshots of the (numerical) solutions over a given initial time period.  Each solution trajectory is based on perturbed initial conditions for \cref{eq: conservation_law}.

\begin{rem} \label{rem:observations}
Our underlying assumption is that realizable training data come from noisy or sparse observations of the conserved state variables of the true underlying physical system. For practical reasons, in this investigation the observations are obtained by numerically simulating the true PDE from perturbed initial conditions. In particular, since we have no oracle knowledge of how the solution may evolve, we do not include initial conditions that may more readily yield a solution space with piecewise smooth profiles.  As was noted in \cite{chencfn}, this is a significant departure from other algorithms that include a large class of initial conditions to simulate  trajectories of conserved state variables \cite{chenPDE}. \end{rem}

We stress that the flux functions in the system corresponding to each conserved quantity are  {\em not} known apriori. Our goal is to learn the dynamics of the conserved variables of the  hyperbolic system and to predict their evolution after this initial time period. In what follows we use the notations adopted from \cite{chencfn} to set up the problem as well as to provide details of how the observations are obtained. For ease of presentation we describe our new methodology for the 1D scalar conservation law, although our numerical experiments in \Cref{sec: numerical} include 1D systems and the 2D Burgers equation. 

For simplicity in both training and testing we use a uniform time step $\Delta t$. We assume that for some chosen $L$ we have available to us $N_{traj}$ observations  over the time period
\begin{equation}
    \label{eq:totalperiod}
    \mathcal{D}_{train}= [0,L\Delta t],
\end{equation} 
which for this investigation will be numerically simulated solutions for our PDEs (see \Cref{sec: pdeexamples}) for $N_{traj}$ initial conditions.  We then break up \cref{eq:totalperiod} into smaller training domains.  Specifically, for the $k$th trajectory where $k = 1,\cdots, N_{traj}$, we assume that training data are available at a set of discrete time instances $\{t_l^{(k)}\}_{l=0}^{L_{train}}$ for a specified $L_{train}$ corresponding to the training period 
\begin{equation}
    \label{eq:trainingperiod}
    \mathcal{D}^{(k)}_{train}= [t_0^{(k)},t^{(k)}_{L_{train}}].
\end{equation} 
Here the initial time $t_0^{(k)} \in [0, (L-L_{train})\Delta t]$, and $t^{(k)}_{L_{train}} = t_0^{(k)} + L_{train}\Delta t$.  We write the conserved state variable information as
\begin{equation}
\label{eq: trajectory data}
\bu^{(k)}(t_l) \in {\mathbb R}^{n_{train}}, \ l=0, \cdots, L_{train}, \ k=1, \cdots, N_{traj},
\end{equation}
where $n_{train}$ denotes the number of spatial grid points in \cref{eq:spatialgrid}.\footnote{We omit the superscript $(k)$ on $t_l$ in \cref{eq: trajectory data} to avoid cumbersome notation.}
 We will construct a network to learn the dynamics and predict solutions $\bu(t_l) \in \mathbb{R}^{n_{test}}$ for $l>L$ and some user prescribed $n_{test}$, which may be different from $n_{train}$. For example, if $n_{train} \ll n_{test}$, we know that our method might be useful for compression purposes.  In \Cref{subsec: coarse} we demonstrate our method's consistency with respect to different choices of  \(n_{train}\).  
\subsection{Conservative flux form network}
\label{sec:conservativeformnetwork}
The conservative flux form network proposed in \cite{chencfn} updates the cell average $\wb{u}_j(t_l)$ using \cref{eq: semi-discrete conservative scheme}. The flux term $f_{j + 1/2}(t_l)$ is replaced by the so-called neural flux
\begin{equation}
\label{eq: neural flux}
f_{j + 1/2}^{NN}(t_l) = \cF^{\btheta} \left( \wb{u}_{j - p}(t_l), \cdots , \wb{u}_j(t_l), \cdots, \wb{u}_{j + q}(t_l) \right) := \cF^{\btheta}_{p,q} ( \wb{u}_j(t_l) ),
\end{equation}
where $\cF^{\btheta}$ is a fully connected feed-forward neural network parameterized by several trainable variables collectively denoted by \(\btheta\), and the inputs $\wb{u}_{j - p}(t_l),\cdots, \wb{u}_{j + q}(t_l)$ are neighboring cell averages centered at $x=x_{j-p}, \cdots, x_{j+q}$, respectively. In this case the fully connected neural network is parameterized by trainable variables having $M$ layers for which the input is given by some $\mathbf{v}\in \mathbb{R}^{p+q+1}$. The network parameters $\btheta$ in \cref{eq: neural flux} include the weight matrices $\mathbf{W}^{(m)} \in \mathbb{R}^{d_m \times d_{m-1}}$ for each layer $m = 1, \dots, M$, where \(d_{0} = p + q + 1\) and \(d_{M}\) is the number of unknowns, as well as the bias vectors $\mathbf{b}^{(m)} \in \mathbb{R}^{d_m}$ for each layer $m = 1, \dots, M-1$. For each $m=1, \cdots, M-1$, we then compute
\begin{equation}\label{eq:activation}
\mathbf{h}^{(m)} = \sigma\left(\mathbf{W}^{(m)} \mathbf{h}^{(m-1)} + \mathbf{b}^{(m)}\right), 
\end{equation}
where $\mathbf{h}^{(0)} = \mathbf{v}$ and $\sigma(\cdot)$ is the activation function  to be specified in \Cref{sec:networkdetails}. Finally, the output $\mathbf{y} \in \bbR^{d_M}$ is obtained by
\[
\mathbf{y} = \mathbf{W}^{(M)} \mathbf{h}^{(M-1)}.
\]

As demonstrated in \cite{chencfn},  the conservative flux form neural network helps to retain important properties in the predicted solution, including  the correct shock propagation speed and conservation.  It was not designed to be entropy-stable, however.

\subsection{Entropy-stable network}
\label{sec:entropystablenetwork}

By incorporating entropy-stability into the conservative flux form network framework, we are able to enhance the quality of the predicted solution.  In this regard, since the KT scheme (see \Cref{sec:KTscheme}) is non-oscillatory, and furthermore TVD in the 1D scalar case, it  effectively ``denoises'' the numerical solution.  This means it is more likely to yield predictions that maintain important physical properties of the solution, including being entropy-stable. As will be observed in \Cref{sec: numerical}, such properties are crucial in low signal-to-noise ratio (SNR) observation environments.  Specifically, the modLW scheme (see \Cref{sec:modLWscheme}) {\em does allow} oscillatory behavior in it solutions, which evidently leads to less accurate predictions.  Finally, we note that our approach is not limited to using the KT scheme.  Other entropy-stable schemes may be similarly incorporated, although may require different specifications for their corresponding hyper-parameters.

\subsubsection{KT-enhanced CFN}\label{sec:KTenhancedCFN} Recall that the KT scheme has conservation form \cref{eq: Kurganov Tadmor Scheme} with numerical flux \cref{eq: KT flux}. As an extension to \cref{eq: neural flux}, we incorporate  \cref{eq: KT flux} to generate the {\em KT-enhanced CFN}, which we write as 
\begin{equation}
H^{NN}_{j+1 / 2}(t):=\frac{\cF^{\btheta}\left(u_{j+1 / 2}^{+}(t)\right)+\cF^{\btheta}\left(u_{j+1 / 2}^{-}(t)\right)}{2}-\frac{a^{NN}_{j+1 / 2}(t)}{2}\left[u_{j+1 / 2}^{+}(t)-u_{j+1 / 2}^{-}(t)\right],
\label{eq: Numerical flux NN}
\end{equation}
where 
\begin{equation}
a^{NN}_{j+1 / 2}(t):=\max \left\{\rho\left(\frac{\partial \cF^{\btheta}}{\partial u}\left(u_{j+1 / 2}^{+}(t)\right)\right), \rho\left(\frac{\partial \cF^{\btheta}}{\partial u}\left(u_{j+1 / 2}^{-}(t)\right)\right)\right\} 
\label{eq: Maximum wave speed NN}
\end{equation}
is the approximate maximum wave speed. As in \cref{eq: neural flux}, $\cF^\btheta$ is a fully connected feed-forward neural network operator characterized by trainable variables collectively denoted by $\btheta$.

Determining the spectral radius of the Jacobian matrix, $\ds \frac{\partial \cF^{\btheta}}{\partial u}$ in \cref{eq: Maximum wave speed NN}, is inherently problem dependent since no fast algorithm exists to compute the spectral radius of a general matrix. In this regard we note that  the spectral radius is simply the corresponding absolute value for the scalar problem, and that there is a closed form formula for systems of two state variables. As there is no explicit form or fast algorithm for spectral radius computation in  systems of more than two state variables, here we approximate the local maximum wave speed $a^{NN}_{j+1 / 2}$ in \cref{eq: Maximum wave speed NN} via  neural network approximation, noting that $\ds \frac{\partial \cF^{\btheta}}{\partial u}$ may be accurately computed using automatic differentiation. 
 
Specifically, we apply a new positive output network to learn the metric for determining the wave speed, resulting in
\begin{equation}
\alpha^{NN}_{j+1 / 2} := \max \left\{\left|\rho_w\left(\frac{\partial \cF^{\btheta}}{\partial u}\left(u_{j+1 / 2}^{+}(t)\right)\right)\right|, \left|\rho_w\left(\frac{\partial \cF^{\btheta}}{\partial u}\left(u_{j+1 / 2}^{-}(t)\right)\right)\right|\right\}.
\label{eq: rho_W}
\end{equation}
Here \(\rho_w(\cdot)\) is approximated by a fully connected neural network whose inputs are all of the entries of the Jacobian matrix.  The parameters are collectively denoted by \(w\).

\subsubsection{LW-enhanced CFN} \label{sec:LWenhancedCFN}
The {\em LW-enhanced CFN} can be obtained by similarly replacing \(f(\cdot)\) in \cref{eq: LWS} by \(\mathcal{F}^{\btheta}(\cdot)\), yielding

\begin{equation}
\begin{aligned}
\wb{u}_{j}^{l+1} = \,& \wb{u}_{j}^{l} 
- \frac{1}{2}\frac{\Delta t}{\Delta x} 
\left[ \mathcal{F}^{\btheta}(\wb{u}_{j+1}^{l}) - \mathcal{F}^{\btheta} ( \wb{u}_{j-1}^{l} )\right]  \\
+ \,& \frac{\Delta t^2}{2\Delta x^2} 
\left[ \frac{\partial \mathcal{F}^{\btheta}}{\partial u}
\left( \frac{\wb{u}_{j}^{l} + \wb{u}_{j+1}^{l}}{2} \right)\left(\mathcal{F}^{\btheta}(\wb{u}_{j+1}^l) - \mathcal{F}^{\btheta}(\wb{u}_{j}^l)\right)\right] \\
- \,&\frac{\Delta t^2}{2\Delta x^2} \left[\frac{\partial \mathcal{F}^{\btheta}}{\partial u}\left( \frac{\wb{u}_{j}^{l} + \wb{u}_{j-1}^{l}}{2} \right)\left( \mathcal{F}^{\btheta}(\wb{u}_{j}^l) - \mathcal{F}^{\btheta}(\wb{u}_{j-1}^l)\right)  \right]. 
\end{aligned}
\label{eq: LWS-cfn}
\end{equation}

\subsubsection{modLW-enhanced CFN} \label{sec:modLWenhancedCFN}
The {\em modLW-enhanced CFN} can be analogously derived by replacing \(f(\cdot)\) in \cref{eq: modLWS} by \(\mathcal{F}^{\btheta}(\cdot)\), yielding
\begin{equation}
\begin{aligned}
\wb{u}_{j}^{l+1} = \,& \wb{u}_{j}^{l} 
- \frac{1}{2}\frac{\Delta t}{\Delta x} 
\left[ \mathcal{F}^{\btheta}(\wb{u}_{j+1}^{l}) - \mathcal{F}^{\btheta} ( \wb{u}_{j-1}^{l} )\right] 
+ \frac{1}{2} \frac{\Delta t^2}{\Delta x^2} 
\Delta_{-}\left[ \frac{\Delta_{+} \mathcal{F}^{\btheta}f ( \wb{u}_{j}^{l})}{\Delta_{+} \wb{u}_{j}^{l}} \Delta_{+}\mathcal{F}^{\btheta} ( \wb{u}_{j}^{l} ) \right]\\
&+ \dfrac{\Delta t}{\Delta x} \Delta_{-} \left[ 
C \gamma\left( \frac{|\Delta_{+} \wb{u}_{j}^{l}|}{(\Delta x)^{\alpha}} \right) \Delta_{+} \frac{\partial \mathcal{F}^{\btheta}}{\partial u}\left( \wb{u}_{j}^{l} \right) \Delta_{+} \wb{u}_{j}^{l} \right]. \\
\end{aligned}
\label{eq: modLWS-cfn}
\end{equation}
As already noted, the allowable oscillations in \cref{eq: modLWS} may lead to instability in \cref{eq: modLWS-cfn}. This will be demonstrated in our numerical experiments. 

\subsection{Boundary Conditions}\label{sec:boundaryconditions} For simplicity we assume that the boundary conditions in \cref{eq: conservation_law} are known, with periodic boundary conditions for the 1D  and 2D  Burgers' equations, and  Dirichlet boundary conditions for shallow water  and Euler's  equations (see \Cref{sec: pdeexamples}.)  In this regard, periodic padding is simple and effective to use in the former, while extra ghost points for padding at the boundary are required in the latter. 
 \begin{rem}\label{rem:boundaryconditions}
     Since error in the initial conditions propagates in time, this approach will not work for other types of (known) boundary conditions if not predefined. Moreover, learning unknown boundary conditions from a finite time interval of observations is in general non-trivial. As the goal of this investigation is to develop entropy-stable conservative neural networks, we leave this task for future research.  
 \end{rem}

\subsection{Time integration}
  With the approximate numerical flux \cref{eq: Numerical flux NN} in hand, we now write the neural net form of the KT conserved form in \cref{eq: Kurganov Tadmor Scheme} as
\begin{equation}
\frac{d}{d t} \wb{u}_j(t)=-\frac{H^{NN}_{j+1/2}(t)-H^{NN}_{j-1 / 2}(t)}{\Delta x}.
\label{eq: KT-enhanced NN}
\end{equation}
We solve \cref{eq: KT-enhanced NN} using the method-of-lines approach, and in our numerical experiments employ the  TVDRK3 method (see \Cref{alg: TVDRK3}). We note that no formal stability analysis exists for explicit time-stepping methods for \cref{eq: KT-enhanced NN}, and we choose TVDRK3 because of its convergence properties for classical numerical conservation laws.\footnote{Recall that the LW and modLW methods are simultaneously temporally and spatially discretized, leading to a one-step time integration (see \Cref{sec:LxWscheme} and \Cref{sec:modLWscheme}).} 

\subsection{The Recurrent Loss Function} \label{sec:recurrentloss} 
We annotate the procedure that comprises  \cref{eq: upm}-\cref{eq: minmod limiter} to \cref{eq: Numerical flux NN}-\cref{eq: rho_W} and the TVDRK3 method (\Cref{alg: TVDRK3}) by \(\mathcal{N}_{KT}\), and similarly, \(\mathcal{N}_{LW}\) and \(\mathcal{N}_{mLW}\) will be used to represent the neural net operator for the LW-enhanced CFN method in \cref{eq: LWS-cfn} and the modLW-enhanced CFN method in \cref{eq: modLWS-cfn}, respectively. Given ${\bu}^l \in \mathbb{R}^{n_{train}}$, $l = 0,\dots, L_{train}$, the next step predictions given by these enhanced CFN methods therefore correspondingly have the form
\begin{equation}
    \label{eq:Ndefine}
    \bu_{KT}^{l+1} = \mathcal{N}_{KT}(\bu^{l}), \quad \bu_{LW}^{l+1} = \mathcal{N}_{LW}(\bu^{l}), \quad \bu_{mLW}^{l+1} = \mathcal{N}_{mLW}(\bu^{l}).
\end{equation}
To avoid cumbersome notation we will omit subscripts \(\cdot_{KT}\), \(\cdot_{LW}\), and \(\cdot_{mLW}\)  when the context is clear. Similarly, we will also use $n$ to denote both $n_{train}$ and $n_{test}$ in sequel. More details regarding $\mathcal{N}_{KT}$, \(\mathcal{N}_{LW}\), and $\mathcal{N}_{mLW}$ are provided in \Cref{sec:networkdetails}.

Following \cite{chencfn} we now define the \emph{recurrent loss function} as
\begin{equation}
\label{eq:recurrentloss}
\mathcal{L}\left( \btheta \right) = \dfrac{1}{N_{traj}} \dfrac{1}{L_{train}} \sum_{k=1}^{N_{traj}} \sum_{l=0}^{L_{train}} \left\|\bu_{NN}^{\left( k \right)}\left( t_{l} ;\btheta\right) - \bu^{\left( k \right)}\left( t_{l} \right)\right\|_2^2, 
\end{equation}
for the $N_{traj}$ trajectory datasets in \cref{eq: trajectory data}, where
\begin{equation}\label{eq:neuralupdate}
    \bu_{NN}^{\left( k \right)}\left( t_{l};\btheta \right) = \underbrace{\mathcal{N}\circ\ldots \circ \mathcal{N}}_{l \text{ times }}\left( \bu^{\left( k \right)}\left( t_{0} \right) \right).
\end{equation}
and $\mathcal{N} = \mathcal{N}_{KT}$,  \(\mathcal{N} = \mathcal{N}_{LW}\) or $\mathcal{N} = \mathcal{N}_{mLW}$ according to \cref{eq:Ndefine} as appropriate.

%% file: 4_Design.tex
\section{Experiment design}
\label{sec:experimentdesign}
We now provide details related to data collection and training used for our numerical experiments in \Cref{sec: numerical}.

\subsection{Data collection}
\label{sec:datacollect}
We test our method and validate our results on the 1D and 2D Burgers' equation with periodic boundary conditions,  the 1D shallow water equations, and the 1D Shu-Osher problem \cite{SHU1988439}, which is a system of Euler equations for gas dynamics. Each represents a classical problem that is often considered in the development of numerical conservation laws \cite{Hesthaven,LeVeque92,Tadmor1987}. The numerical solutions simulated by PyClaw \cite{clawpack,pyclaw} will be used  both to generate synthetic training data as well as to compute reference solutions to evaluate our results.  We reiterate that while we use knowledge of the true flux to generate the  training data, for the neural flux approximation we assume only that the underlying PDE is a conservation law, that is, we {\em do not} use any other flux information.

The number of spatial grid points in the training data, $n_{train}$, is first chosen so that the numerical solutions are well-resolved.  The value of $n_{train}$ is decreased in \Cref{subsec: coarse} to see how our method performs when instead given   coarsely-sampled observations.  We fix the time step \(\Delta t\) for each PDE to satisfy the Courant–Friedrichs–Lewy (CFL) condition in the finest resolution (largest $n_{train}$) case.  In so doing we limit the number of parameters in the problem and ensure meaningful comparisons between results. We note that this is not an inherent limitation of our method, and that using variable-sized time steps could improve computational efficiency.

\subsection{Network and training details}
\label{sec:networkdetails}
 The neural flux operator, $\mathcal{N}$ in \cref{eq:neuralupdate}, is constructed using a fully connected neural network, with the input and output size equivalent to the number of state variables in each problem. Following what was done in \cite{chencfn}, our experiments  utilize five hidden layers with 64 hidden neurons at each layer.  We use the silu activation function \(\text{silu}(x) = \text{Sigmoid}(x)\cdot\text{ReLU}(x)\) in \cref{eq:activation} since it is differentiable and therefore appropriate for computing the Jacobian  in \cref{eq: rho_W}.  
The network used to approximate the spectral radius in \cref{eq: rho_W} is a 2-hidden-layer, 64-hidden-neurons-per-layer neural network with activation function \(\text{ReLU}(x)\) (since  second-order differentiability  is not required in this case). 

We employ the commonly used  Adam stochastic optimization method \cite{Adam} to update the neural networks' weights and biases. For consistency we train our numerical models  for 500 epochs with exponential decay learning rate whose initial value is dependent on problem complexity -- $10^{-4}$ for both 1D and 2D Burgers' equation case and $2\times10^{-3}$ for both the shallow water and Euler's equation examples. We fix the batch size at 10 for all 1D experiments.  In consideration of GPU memory limits the batch size for the 2D Burgers' equation experiment is 1. While the number of training trajectories, \(N_{traj}\), varies by type of experiment, we fix the number of validation cases (which determines whether to save the updated parameters at each epoch) at 40. All implementations are based on JAX \cite{jax}, a high-performance numerical computing library that provides automatic differentiation and GPU/TPU acceleration for Python. Finally, to ensure robustness, we note that none of these parameters were fine-tuned.\footnote{The complete code  is  available upon request for reproducibility purposes.} 

\subsection{Conservation and entropy metric}\label{sec:conservationnmetric}
In addition to evaluating our new entropy-stable CFN predictions by comparing them to corresponding reference solutions occurring beyond the training domain \cref{eq:trainingperiod}, we also seek to demonstrate the preservation of important  conservation law properties, namely the conservation of each conserved quantity and entropy stability. In this regard we explicitly define the discrete conserved quantity remainder at time $t_l$ as 
\begin{equation}
\label{eq:conservemetric}
C(\bu(t_l)) := \left| \sum_{j=1}^{n-1} \left( \bar{u}_j(t_l) - \bar{u}_j(t_0) \right)\Delta x - \sum_{s=1}^{l}\left( F_a^{s-1} -F_b^{s-1} \right)\Delta t \right|,
\end{equation}
where $\bu(t) = (\bar{u}_0(t), \ldots, \bar{u}_{n}(t))^T \in \mathbb{R}^n$ is the prediction at time \(t\), and flux terms \(F_a^{s-1}\) and \(F_b^{s-1}\) are the calculated flux operators at each respective boundary defined by
\begin{equation}
\label{eq:fluxterm}
        F_a^{s-1} = \frac{1}{\Delta t}\int_{t_{s-1}}^{t_{s}} f(u(a,t))dt,\quad
        F_b^{s-1} = \frac{1}{\Delta t}\int_{t_{s-1}}^{t_{s}} f(u(b,t))dt.
\end{equation}
Since the fluxes are unknown, we use $F^\theta$ as defined in \eqref{eq: neural flux} to calculate \cref{eq:fluxterm}. 

As for the entropy, all of our prototypical examples have explicit forms for their (typically used) entropy pairs (see e.g.~\cite{LeVeque92}). Entropy evolution for each example is therefore evaluated by simply substituting the predictions into the corresponding entropy functions at each time step. Specifically, the discrete entropy remainder term analogous to \cref{eq: entropy inequality semi} is given by
\begin{equation}
\mathcal{J}(\bu (t_l)) := \left(  \sum_{j=1}^{n-1} \left( \bar{U}_j(t_l) - \bar{U}_j(t_0) \right)\Delta x - \sum_{s=1}^{l}\left( \bar{F}_a^{s-1} -\bar{F}_b^{s-1} \right)\Delta t\right), 
    \label{eq: discrete entropy}
\end{equation}
where \(\bar{U}\) is the entropy function and \(\bar{F}_a, \bar{F}_b\) are the evaluations of the corresponding entropy flux at each respective boundary. We say we have obtained an {\em entropy-stable} network operator if $\mathcal{J}(\bu (t_l)) \leq 0$.

Finally, we also measure the relative \(\ell^2\) prediction error, which we define as
\begin{equation}
    \label{eq:relative l2 error}
    \mathcal{R}(\bu(t_l), \bu_{true}(t_l)) := \frac{\left\| \bu(t_l) - \bu_{true}(t_l)\right\|_2}{\left\| \bu_{true}(t_l) \right\|_2}.
\end{equation}
Here  $\bu_{true}(t)$ is the true solution at time \(t\) and we compute the discrete $\ell^2$ norm as \(\left\| \bu \right\|_2 = \sqrt{{\Delta x}\sum_{j} \bar{u}_j^2}\). 

A few remarks are in order:

\begin{remark}
    \label{rem:metrics} We  emphasize that the metric \cref{eq:conservemetric} uses the {\em learned} flux and not the true flux known in each of our experiments.  This is because we want to determine whether or not the enhanced CFN method maintains conservation for the prediction, as opposed to the model itself.  However the enhanced CFN does not learn each corresponding entropy pair, so in this case we use what is learned to compute the true entropy pair for each underlying model. This is not incorporated into the prediction, and is used solely for the purpose of evaluating the enhanced CFN's ability to satisfy entropy stability.
\end{remark}

\begin{remark}
\label{rem:relative l2 error}
The relative \(\ell^2\) prediction error, as defined in \cref{eq:relative l2 error}, accounts solely for the {\em global} behavior of the solution.  This may yield misleading outcomes, especially in cases where the solutions admit discontinuous or oscillatory profiles. A more thorough discussion is reserved for the numerical experiments conducted in \Cref{sec: numerical}. 
\end{remark}

%% file: 5_Numerical.tex
\section{Numerical Experiments}
\label{sec: numerical}

We now provide a series of numerical experiments to demonstrate the benefits of the entropy stable CFN.  In particular we will observe the critical importance of the non-oscillatory property in the KT algorithm \cite{KTscheme2000} for long term stability, as our examples will show that the KT-enhanced CFN method is more robust than both the LW-enhanced CFN method and modLW-enhanced CFN method. All of our entropy stable neural networks yield better accuracy and stability when compared to the original CFN method in \cite{chencfn}, however. Since they are often used for validating the efficacy of numerical methods for conservation laws, we will consider Burgers' equation (in 1D and 2D), the 1D system of shallow water equations, and the 1D system of Euler's equations as characterized in the Shu-Osher shock tube problem \cite{SHU1988439}. We review each of these classical problems in \Cref{sec: pdeexamples}.

In  \Cref{subsec: noisy} we conduct numerical experiments for each example where we assume we have available dense observations that are corrupted by noise of varying amounts ranging from $0$ to $100\%$. We then analyze the temporal evolution of the relative $\ell^2$ prediction errors for each conserved variable.
In \Cref{subsec: coarse} we investigate our method's robustness and consistency by training on coarser spatial resolutions and then testing on fully resolved data. Finally,  the sensitivity analysis in \Cref{subsec: sensitivity} determines the impact of alternative specifications on our approach.

\subsection{Prototype Conservation Laws Used in Numerical Experiments}
\label{sec: pdeexamples}
\subsubsection{1D Burgers' Equation}
\label{sub:1DBurgers}
We first consider scalar Burgers' equation
\begin{equation}
    \dfrac{\partial u}{\partial t} + \dfrac{\partial}{\partial x} \left( \dfrac{u^2}{2} \right) = 0, \quad x \in \left[ 0, 2\pi \right], \quad t > 0,
    \label{eq: Burgers' equation}
\end{equation}
with periodic boundary conditions $u(0,t)=u(2 \pi, t)$.  The initial conditions for the training data are given by 
\begin{equation}
    u\left( x, 0 \right) = \alpha \sin\left( x \right) + \beta,
    \label{eq: Burgers' initial condition}
\end{equation}
where $\alpha$ and $\beta$ are sampled uniformly from intervals $\left[ .75, 1.25 \right]$ and $\left[ -.25, .25 \right]$, respectively.
As discussed in \Cref{sec:datacollect}, each of the $N_{traj}$ training datasets is generated  using the  PyClaw package \cite{clawpack, pyclaw}, for which we fix time step \(\Delta t = .005\). 
We also set $N_{traj}=200$, with the total training time period \cref{eq:totalperiod} defined by $L = 20$, i.e.~$[0,.1]$, and each temporal training domain \cref{eq:trainingperiod} defined by $L_{train} = L =20$. That is, in this first experiment we do not break up the total training period into smaller temporal training domains.
  We emphasize that only {\em smooth} data profiles, that is, before the shock forms, are fed into the training procedure.

For the testing data we fix $\alpha$ and $\beta$ in \cref{eq: Burgers' initial condition} and consider the specific initial conditions
 \(u_{test}\left( x,0 \right) = 1.05609\sin x + 0.1997.\)
The reference solution is then generated up to $T = 600 \Delta t = 3$ using PyClaw  with the same time step as the training data and with grid size \(n_{test} = 512\). 
As is standard, we define the entropy pair $(U(u),F(u))$ for \cref{eq:entropy inequality} to be used in \cref{eq: discrete entropy} as $(\frac{u^2}{2}, \frac{u^3}{6})$. 

\subsubsection{Shallow Water Equation} \label{sub:shallowwater} The $1D$ system of shallow water equation considered is given by
\begin{equation}
    \begin{aligned}
        \dfrac{\partial h}{\partial t} + \dfrac{\partial}{\partial x} \left( hu \right) &= 0, \\
        \dfrac{\partial}{\partial t} \left( hu \right) + \dfrac{\partial}{\partial x} \left( hu^2 + \dfrac{1}{2}gh^2 \right) &= 0.
    \end{aligned}
    \label{eq: shallow water equation}
\end{equation}
 We assume Dirichlet boundary conditions and initial conditions given by 
\[
    h\left( x, 0 \right) = \left\{
    \begin{aligned}
     &h_l + \xi_{l}, \quad x < x_0 + \xi_{x}   \\
     &h_r + \xi_{r}, \quad x \geq x_0 + \xi_{x} \\
     \end{aligned}
    \right.,
    u\left( x, 0 \right) = \left\{
    \begin{aligned}
     &u_l + \xi_{ul}, \quad x < x_0 + \xi_{x}   \\
     &u_r + \xi_{ur}, \quad x \geq x_0 + \xi_{x} \\
     \end{aligned}
    \right..
\]
Here \(h_l = 3.5, h_r=1.0, u_l=u_r=x_0=0\), with \(\xi_{l}, \xi_{r}\in \mathcal{U}\left[ -.2, .2 \right]\) and \(\xi_{ul}, \xi_{ur}, \xi_{x} \in \mathcal{U}\left[ -.1,.1 \right]\) sampled uniformly.   The spatial domain considered is $(-5,5)$.

Once again the training data and reference solutions are computed using PyClaw \cite{clawpack,pyclaw}  (HLLE Riemann Solver). We generate $N_{traj} = 200$ training trajectories and use time step \(\Delta t = .005\) to satisfy the CFL condition for any of the spatial grid values of  $\Delta x$ used in our experiments. The total training time period \cref{eq:totalperiod} is defined by $L = 20$, i.e.~$[0,.1]$, with each temporal training domain  \cref{eq:trainingperiod} defined by $L_{train} = L = 20$.

The testing data is generated by setting \(h_l = 3.5691196, h_r = 1.178673, u_l = -.064667,\) \(u_r = -.045197, x_0 = .003832\) with grid size $n_{test} = 512$ and $\Delta t = .005$ up to \(T = 200\Delta t = 1\). The entropy function is equivalent to the the total energy of the system, \(U = \frac{1}{2}gh^2 + \frac{1}{2}hu^2\), with corresponding entropy flux $F = \frac{1}{2}hu^3 + gh^2u$.   These values are used for the  entropy metric \cref{eq: discrete entropy}. 

\subsubsection{Euler's Equation}\label{sub:Eulers} The system of 1D Euler's equations is expressed by
\begin{equation}
    \begin{aligned}
        \rho_{t} + \left( \rho u \right)_{x} &= 0, \\
        \left( \rho u \right)_{t} + \left( \rho u^2 + p \right)_{x} &= 0, \\
        \left( E \right)_{t} + \left( u (E + p) \right)_{x} &= 0.
    \end{aligned}
    \label{eq: Euler's equation}
\end{equation}
Here again we consider the spatial domain $(-5,5)$
and assume Dirichlet boundary conditions, with initial conditions given by
\[
\begin{aligned}
& \rho(x, 0)= \begin{cases}\rho_l, & \text { if } x \leq x_0, \\
1+\varepsilon \sin (5 x), & \text { if } x_0<x \leq x_1, \quad u(x, 0)= \begin{cases}u_l, & \text { if } x \leq x_0, \\
0, & \text { otherwise, }\end{cases} \\
1+\varepsilon \sin (5 x) e^{-\left(x-x_1\right)^4} & \text { otherwise, }\end{cases} \\
& p(x, 0)=\left\{
\begin{array}{ll}
p_l, & \text { if } x \leq x_0, \\
p_r, & \text { otherwise, }
\end{array} \quad E(x, 0)=\frac{p_0}{\gamma-1}+\frac{1}{2} \rho(x, 0) u(x, 0)^2 . \right.
\end{aligned}
\]
The parameters are sampled uniformly from 
$$
\begin{aligned}
& \rho_l \sim \mathcal{U}\left[\hat{\rho}_l(1-\epsilon), \quad \hat{\rho}_l(1+\epsilon)\right], \quad \hat{\rho}_l=3.857135, \\
& \varepsilon\; \sim \mathcal{U}\left[\hat{\varepsilon}(1-\epsilon), \quad \hat{\varepsilon}(1+\epsilon)\right], \quad \hat{\varepsilon}= .2, \\
& p_l \sim \mathcal{U}\left[\hat{p}_l(1-\epsilon), \quad \hat{p}_l(1+\epsilon)\right], \quad \hat{p}_l=10.33333, \\
& p_r \sim \mathcal{U}\left[\hat{p}_r(1-\epsilon), \quad \hat{p}_r(1+\epsilon)\right], \quad \hat{p}_r=1, \\
& u_l \sim \mathcal{U}\left[\hat{u}_l(1-\epsilon), \quad \hat{u}_l(1+\epsilon)\right], \quad \hat{u}_l=2.62936, \\
& x_0 \sim \mathcal{U}\left[\hat{x}_0(1-\epsilon), \quad \hat{x}_0(1+\epsilon)\right], \quad \hat{x}_0=-4, \\
&
\end{aligned}
$$
with $\epsilon=.1, x_1=3.29867$ and $\gamma=1.4$. Noting that $\hat{\rho}$, $\hat{p}_l$, $\hat{u}$ are the same values as those used in the CLAWPACK Shu-Osher example, the $N_{traj} = 300$ training trajectories are generated by employing PyClaw (HLLE Riemann Solver) with  \(\Delta t = .002\) to solve the system. The total training period $\mathcal{D}_{train}$ in \cref{eq:totalperiod} is defined by  \(L = 300\), yielding $[0,.6]$, while the training period $\mathcal{D}_{train}^{(k)}$, $k = 1,\dots, N_{traj}$, in \cref{eq:trainingperiod} is defined by $L_{train} = 20$ with $t_0^{(k)} = t_l$,  $l \in \{0,1,\ldots, L - L_{train}\}$. 

The testing data are generated by fixing the parameters \(\rho_l = \hat{\rho}_l, \varepsilon = \hat{\varepsilon}, p_l = \hat{p}_l, p_r = \hat{p}_r, u_l = \hat{u}_l,\) and \(x_0 = \hat{x}_0\). The time step is again $\Delta t = .002$ with simulations conducted up to \(T = 800\Delta t = 1.6\). The spatial grid size is chosen so that \(n_{test} = 512\). 

The entropy function is given by
\(
U(\rho, u, p) = -\rho S,
\)
where \(S = log(p/\rho^\gamma)\) represents the specific entropy and \(p\) is the pressure obtained by \(p = (\gamma - 1)E - (\gamma - 1)\frac{\rho u^2}{2}\), while the 
 associated entropy flux is \(
F(\rho, u, E) = -\rho u S
\). The entropy pair $(U,F)$ is used for the entropy metric \cref{eq: discrete entropy}. 

\subsubsection{2D Burgers' Equation}\label{sub: 2D Burgers} The 2D periodic Burgers' equation is given by 
\begin{equation}
    \dfrac{\partial u}{\partial t} + \dfrac{\partial}{\partial x} \left( \dfrac{u^2}{2} \right) + \dfrac{\partial}{\partial y} \left( \dfrac{u^2}{2} \right) = 0, \quad (x, y) \in \left[ 0, 1 \right]\times\left[ 0, 1 \right], \quad t > 0, 
    \label{eq: Burgers' equation 2d}
\end{equation}
subject to  boundary conditions $u(0,y, t)=u(1, y, t)$ and $u(x, 0, t)=u(x, 1, t)$. We specify the initial conditions as
\[
    u\left( x, y, 0 \right) = \alpha \sin\left( 2\pi x +x_0\right)\cos\left( 2\pi y +y_0\right) + \beta,
    \label{eq: Burgers' initial condition 2d}
\]
where
\(\alpha \sim {\mathcal U}[.75,1.25],\beta \sim {\mathcal U}[-.25,.25],  x_0 \sim {\mathcal U}[.5,1.5], \text{ and} y_0 \sim {\mathcal U}[-.5, .5].\)  
The $N_{traj} = 5$ training trajectories and corresponding reference solutions are computed using PyClaw \cite{clawpack,pyclaw} with the HLLE Riemann Solver.  To ensure compliance with the CFL condition for the varying choices of spatial grid $\Delta x$ used in our experiments, we fix our time step as \(\Delta t = .0005\). The total training domain $\mathcal{D}_{train} = [0,.01]$ is  determined by  \cref{eq:totalperiod} with $L = 20$.  We also choose $L_{train} = L = 20$ for each  temporal training domain  $\mathcal{D}_{train}^{(k)}$ in 
 \cref{eq:trainingperiod}, $k = 1,\dots, N_{traj}$. 

We set the parameters for the testing data as
\(x_0 = 1.032833, y_0 = .034137, \alpha = 1.004777, \text{ and }  \beta = .106782.\)
The grid size is chosen so that $n_{test} = 200\times 200$, or equivalently $\Delta x = \Delta y = \frac{1}{200}$. We extend the testing time domain  to \( T = 1600\Delta t = .8\). The entropy pair $(U,F)$ is the same as the 1D case in \Cref{sub:1DBurgers}.

\subsection{Observations with varying noise coefficients}
\label{subsec: noisy}
We proceed by examining the case for which the training data are contaminated by varying amounts of noise as determined by the noise coefficient $\eta \in [0,1]$ in each respective initial condition, \cref{eq: noise burgers}, \cref{eq:noise_shallow}, and \cref{eq:noise_euler}, for the models discussed in \Cref{sec: pdeexamples}.  Initial conditions for the 2D Burgers' equation is analogous to its 1D formulation.  The grid size for training data is fixed at $n_{train} = 512$ for each 1D problem, while we employ 
\(n_{train} = 200\times200\) in the 2D Burgers' equation to accommodate GPU memory limitations.  We show results for the KT-enhanced CFN network \cref{eq: Numerical flux NN}.  Later in \Cref{subsec: sensitivity} we discuss the implications of using the corresponding LW and modLW schemes.

\subsubsection{1D Burgers' Equation} 
\label{subsubsec: 1D Burgers Noise}
Using the specifications defined in \Cref{sub:1DBurgers}, we assume the perturbed training data in domain $\mathcal{D}^{(k)}_{train}$ defined in \cref{eq:trainingperiod} for $k = 1,\dots,N_{traj}$ contain noise of various levels as determined by
\begin{equation}
\tilde{u}^{(k)}\left( x_i, t_l \right) = u^{(k)}\left( x_i, t_l \right) + \eta \overline{|u\left( x,t \right)|} \xi_{i,l},
\label{eq: noise burgers}
\end{equation}
 where $\xi_{i,l} \sim \mathcal{N}\left( 0, 1 \right)$, $i = 1,\dots, n_{train}$, and  $l = 0,\dots, L_{train}$.\footnote{Recall in this example we choose $L_{train} = L$ in \cref{eq:totalperiod}.} 
 Here \(\overline{|u\left( x,t \right)|}\) represents the mean absolute value of exact solution \(u\left( x,t \right)\) over the entire dataset. The noise coefficient \(\eta \in [0,1]\) describes the noise intensity applied to the data. 
\begin{figure}[h!]
    \centering
    \includegraphics[width=0.9\linewidth]{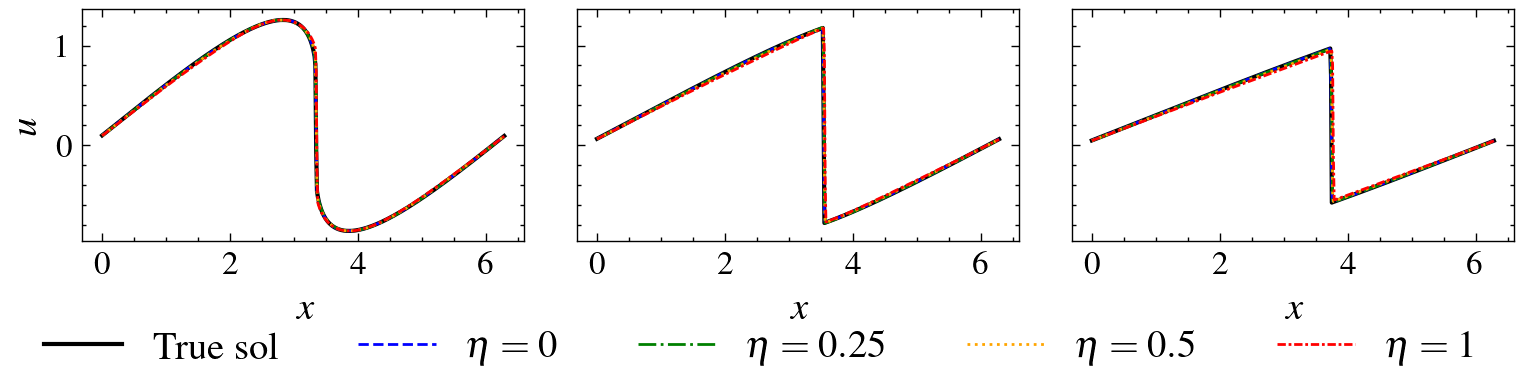}
    \caption{Comparison of the reference solution (black solid line) to 1D Burgers' equation with predictions from the KT-enhanced CFN solutions at (left) \(t = 1\) (middle) \(t = 2\) (right) \(t = 3\) for $\eta = 0,.25,.5,1$ in \cref{eq: noise burgers}.}
    \label{fig: Burgers Noise}
\end{figure}

\Cref{fig: Burgers Noise} shows a temporal sequence of the KT-enhanced CFN solution described in \Cref{sec:KTenhancedCFN} for the 1D Burgers' equation given training data specified by \cref{eq: noise burgers} with  $\eta = 0,.25, .5$ and $1$. Observe that the KT-enhanced CFN solution accurately captures the shock formation which occurs at \( t = 1.0\), {\em even though} the training data are only  collected in $\mathcal{D}_{train} = [0, .1]$, that is prior to the shock formation and while the solution is still smooth.  This means that the method does not require some oracle knowledge of what is needed in the eventual solution space. Moreover, the shock characteristics are maintained for long term prediction, even for $\eta = 1$. 

\begin{figure}[h!]
    \centering
    \includegraphics[width=0.9\linewidth]{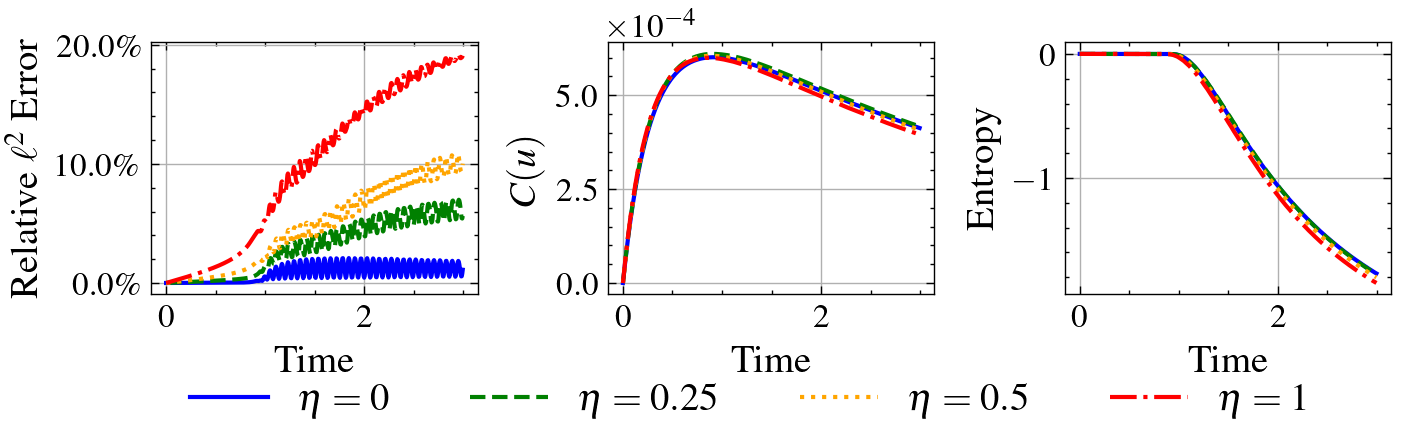}
    \caption{(left) Relative $\ell^2$ prediction error \cref{eq:relative l2 error}, (middle) Discrete conserved quantity remainder $C(u)$ in \cref{eq:conservemetric}, (right) discrete entropy remainder $\mathcal{J}(u)$ in \cref{eq: discrete entropy} for 1D Burgers' equation, $t \in [0,3]$, and $\eta = 0, .25, .5, 1$ in \cref{eq: noise burgers}. }
    \label{fig: Burgers Noise Error}
\end{figure}

There is, however, measurable impact on the overall error when the training data contain significant levels of noise.  \Cref{fig: Burgers Noise Error} (left) shows the relative \(\ell^2\) prediction error for varying noise coefficients in the time domain $[0,3]$. Moreover noise increasingly impacts the solution  once the shock is formed. 

The robustness of KT-enhanced CFN method with respect to noise is further demonstrated in \Cref{fig: Burgers Noise Error} (middle), which displays the evolution of the discrete conserved quantity remainder $C({u})$ in \cref{eq:conservemetric} for  various noise coefficients $\eta$ in \cref{eq: noise burgers}.  \Cref{fig: Burgers Noise Error} (right) shows the corresponding  evolution of  the discrete entropy remainder $\mathcal{J}({u})$ in \cref{eq: discrete entropy}, demonstrating its non-positivity for each choice of $\eta$. In combination these results demonstrate that the KT-enhanced CFN  yields consistent and robust long-term predictions for the 1D Burgers' equation.

\subsubsection{Shallow water Equation} 
\label{subsubsec: Shallow Water Noise}
\begin{figure}[h!]
    \centering
    \includegraphics[width=1\linewidth]{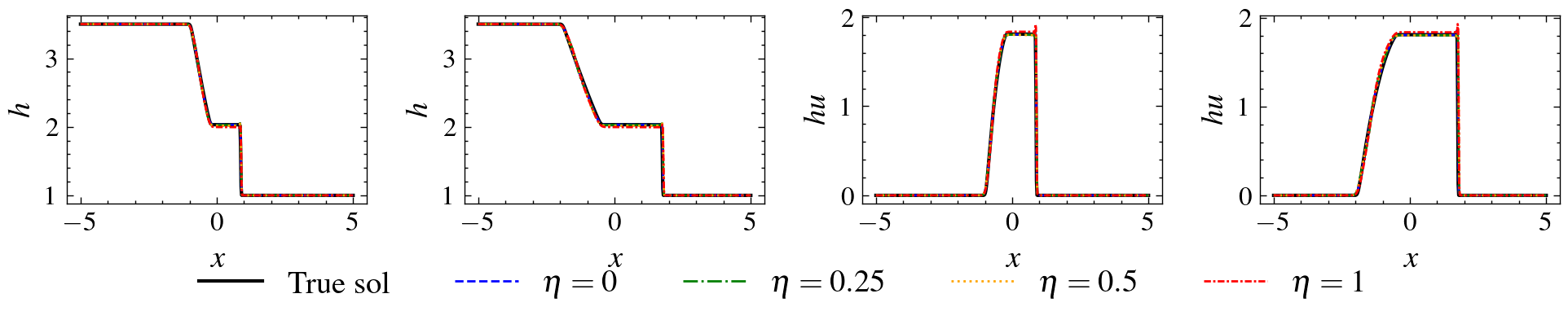}
     \caption{Comparison of the reference solution (black solid line) of height \(h\) and momentum \(hu\) for the shallow water equations with the  KT-enhanced CFN predictions for  \(\eta = 0, .25, .5, 1\) in \cref{eq:noise_shallow}: (left) \(t = .5\) of \(h\), (middle-left) \( t= 1\) of \(h\), (middle-right) \(t = .5\) of \(hu\), (right) \(t = 1\) of \(hu\).}
     \label{fig: Shallow Water Noise h}
\end{figure}
We now examine the effects of noise in the training data for the shallow water equations \cref{eq: shallow water equation}.  As in the 1D Burgers' equation case, zero mean Gaussian noise is added to the training data in the domain $\mathcal{D}^{(k)}_{train}$ defined in \cref{eq:trainingperiod} for $k = 1,\dots,N_{traj}$, yielding
\begin{equation}
    \begin{bmatrix} 
    \tilde{h}\left( x_{i}, t_{l} \right) \\ 
    \tilde{hu}\left( x_{i}, t_{l} \right)
    \end{bmatrix} 
     = 
    \begin{bmatrix}
    h\left( x_{i}, t_{l} \right) \\
    hu\left( x_{i}, t_{l} \right)
    \end{bmatrix}
     + \eta\overline{\bm{u}}
       \xi_{i,l}, \quad \xi_{i,l} \sim \mathcal{N}\left( \begin{bmatrix} 
       0\\
       0
       \end{bmatrix} 
       , \begin{bmatrix} 
                1&0\\
                0&1
       \end{bmatrix} 
       \right),
       \label{eq:noise_shallow}
\end{equation}
where \( i = 1,\dots, n_{train}, l = 0,\dots, L_{train},\)\footnote{Recall in this example we also choose $L_{train} = L$ in \cref{eq:totalperiod}.} and \(\overline{\bm{u}}\) is the mean absolute value of the training data \(\bm{u} = \left[ h, hu \right]^{T}\) over the entire dataset.  The noise intensity coefficient \(\eta \in [0,1]\) is chosen as \(\eta = 0, .25, .5, 1\).

\begin{figure}[h!]
\centering
\includegraphics[width=1.\linewidth]{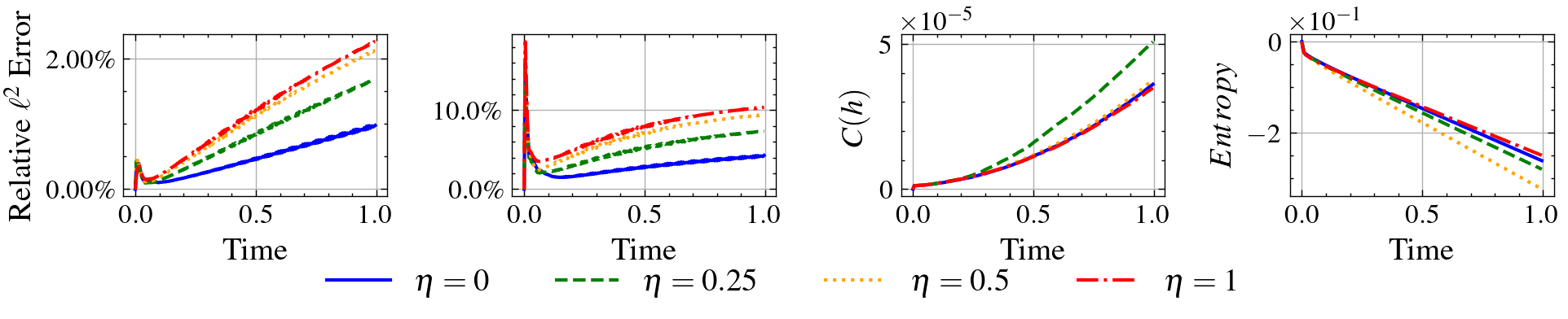}
 \caption{Relative \(\ell^2\) prediction error \cref{eq:relative l2 error} in the shallow water equations: (left) \(h\), (middle-left) \(hu\); and (middle-right) discrete conserved quantity remainder \cref{eq:conservemetric} of height $C(h)$; (right) discrete entropy remainder \(\mathcal{J}([h,hu]^T)\) in \cref{eq: discrete entropy} for the shallow water equations with $\eta = 0,.25,.5, 1$ in \cref{eq:noise_shallow}.}
    \label{fig: Shallow Water Noise Error}
\end{figure}

\Cref{fig: Shallow Water Noise h} displays the height \(h\) and momentum \(hu\) at \(t = .5,\) and \(1\) for each different $\eta$ in \cref{eq:noise_shallow}.  Observe that the KT-enhanced CFN method can accurately predict the shock structure even when  $\eta = 1$.  \Cref{fig: Shallow Water Noise Error} (left and middle-left) depicts the impact of noise over time.  The spike seen in \Cref{fig: Shallow Water Noise Error} (middle-left) is due to $h u(x,0) = 0$ causing a large relative error in the first time step, and demonstrates the sensitivity of the relative \(\ell^2\) prediction error to discontinuity formation  (see \Cref{rem:relative l2 error}). 

\Cref{fig: Shallow Water Noise Error} (middle-right) shows the evolution of the discrete conserved quantity remainder $C(h)$  in \cref{eq:conservemetric}, while the discrete entropy remainder \cref{eq: discrete entropy} is displayed in \Cref{fig: Shallow Water Noise Error} (right).  Observe conservation is preserved up to $\mathcal{O}(10^{-5})$ and that the discrete entropy remainder $\mathcal{J}([h,hu]^T)$ is non-positive. The results for $C(hu)$ are comparable and therefore not displayed.

\subsubsection{Euler's Equation} 
\label{subsubsec: Euler Noise}
We now investigate the effects of noise in the training data for Euler's equations \cref{eq: Euler's equation}.  Once again zero mean Gaussian noise is added to the training data in domain \(\mathcal{D}_{train}^{(k)}\) defined in \cref{eq:trainingperiod} for \(k = 1,\ldots, N_{traj}\), yielding
\begin{equation}
    \begin{bmatrix} 
    \tilde{\rho}\left( x_{i}, t_{l} \right) \\ 
    \tilde{\rho u}\left( x_{i}, t_{l} \right) \\
    \tilde{E}\left(x_{i}, t_{l}\right)
    \end{bmatrix} 
     = 
    \begin{bmatrix}
    \rho\left( x_{i}, t_{l} \right) \\
    \rho u\left( x_{i}, t_{l} \right) \\
    E\left( x_{i}, t_{l} \right)
    \end{bmatrix}
     + \eta\overline{\bm{u}}
       \xi_{i,l}, \quad \xi_{i,l} \sim \mathcal{N}\left( \begin{bmatrix} 
       0\\
       0\\
       0
       \end{bmatrix} 
       , \begin{bmatrix} 
                1&0&0\\
                0&1&0\\
                0&0&1\\
       \end{bmatrix} 
       \right),
       \label{eq:noise_euler}
\end{equation}
where \(i = 1, \ldots, n_{train}\),  \(l = 0, \ldots, L_{train}\),\footnote{Recall in this case \(L_{train} < L\). (See \Cref{sub:Eulers}.)} and \(\overline{\bm{u}}\) is the mean absolute value of the training data \(\bm{u} = \left[ \rho, \rho u, E \right]^{T}\) over the entire dataset up to time \(t = t_{L}\). The noise intensity coefficient  \(\eta \in [0,1]\),  is chosen as \(0, .25, .5,\) and \(1\). 
\begin{figure}[h!]
    \centering
    \includegraphics[width=1.0\linewidth]{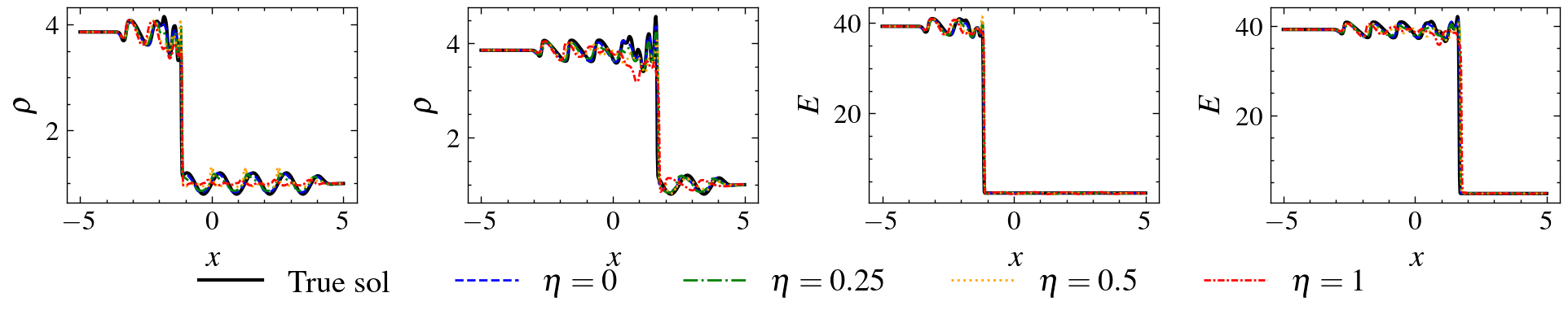}
    \caption{Comparison of the reference solution (black solid line) of density \(\rho\) and energy \(E\) in Euler's equation with the  KT-enhanced CFN predictions with \(\eta = 0, .25, .5, 1\) in \cref{eq:noise_euler}: (left) \(t = .8\) of \(\rho\), (middle-left) \(t = 1.6\) of \(\rho\), (middle-right) \(t =.8\) of \(E\), (right) \(t = 1.6\) of \(E\). }
    \label{fig: Euler Noise rho}
\end{figure}
\begin{figure}[h!]
    \centering
    \includegraphics[width=1.\linewidth]{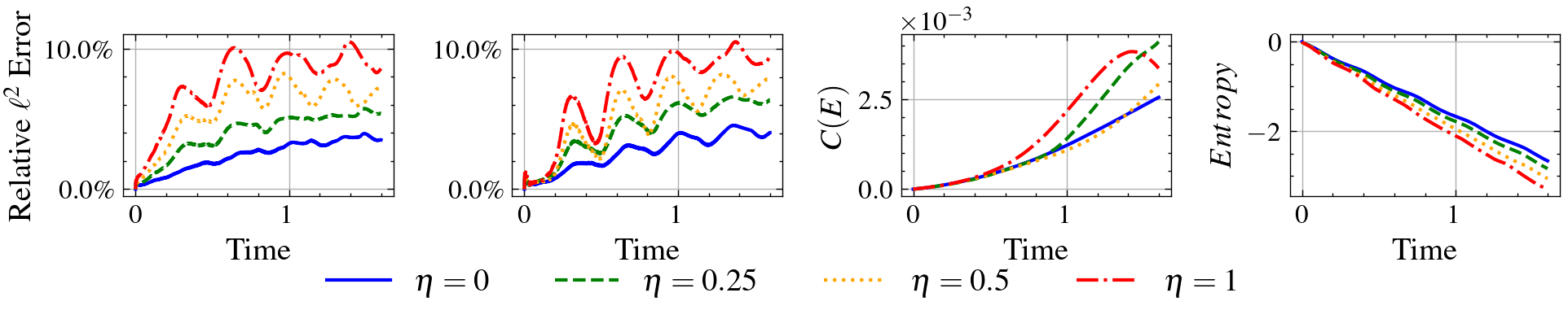}
    \caption{Relative $\ell_2$ prediction error \cref{eq:relative l2 error}: (left) \(\rho\) (middle-left) \(E\), and (middle-right) discrete conserved quantity remainder \(C({E})\) of \cref{eq:conservemetric}, and (right) discrete entropy remainder $\mathcal{J}([\rho, \rho u, E]^T)$ in \cref{eq: discrete entropy} for Euler's equation with \(\eta = 0, .25, .5, 1\) in \cref{eq:noise_euler}.}
    \label{fig: Euler Noise Error}
\end{figure}

\Cref{fig: Euler Noise rho} compares the KT-enhanced CFN predictions of density \(\rho\) and energy \(E\)  for different choices of $\eta$ in \cref{eq:noise_euler}. We omit the comparison of momentum \(\rho u\) due to its similarity to energy \(E\). The improved accuracy is apparent when compared to the original non-entropy stable  CFN \cite[section 5.3, Fig. 20]{chencfn}. The impact of noise is also clearly evident, which is likely due to (i) the additional complexity in the corresponding flux terms of Euler's equation and (ii) the use of the surrogate \cref{eq: rho_W} in place of the true spectral radius in the maximum wave speed approximation \cref{eq: Maximum wave speed NN}, which is not as robust to noise due to this added complexity. Nevertheless, the evolution of relative \(\ell^2\) prediction error shown in \Cref{fig: Euler Noise Error} still demonstrates consistency, even when $\eta = 1$ in \cref{eq:noise_euler}. The robustness and entropy stable property of KT-enhanced CFN with respect to noise is further displayed by the evolution of  the discrete conserved quantity remainders  \(C(E)\), with similar results for \(C(\rho)\)  and  \(C(\rho u)\), and the discrete entropy remainder $\mathcal{J}([\rho, \rho u, E]^T)$ for different noise coefficients. 

\subsubsection{2D Burgers' Equation} 
\label{subsubsec: 2D Burgers Noise}
\begin{figure}[h!]
    \centering
    \includegraphics[width=.9\linewidth]{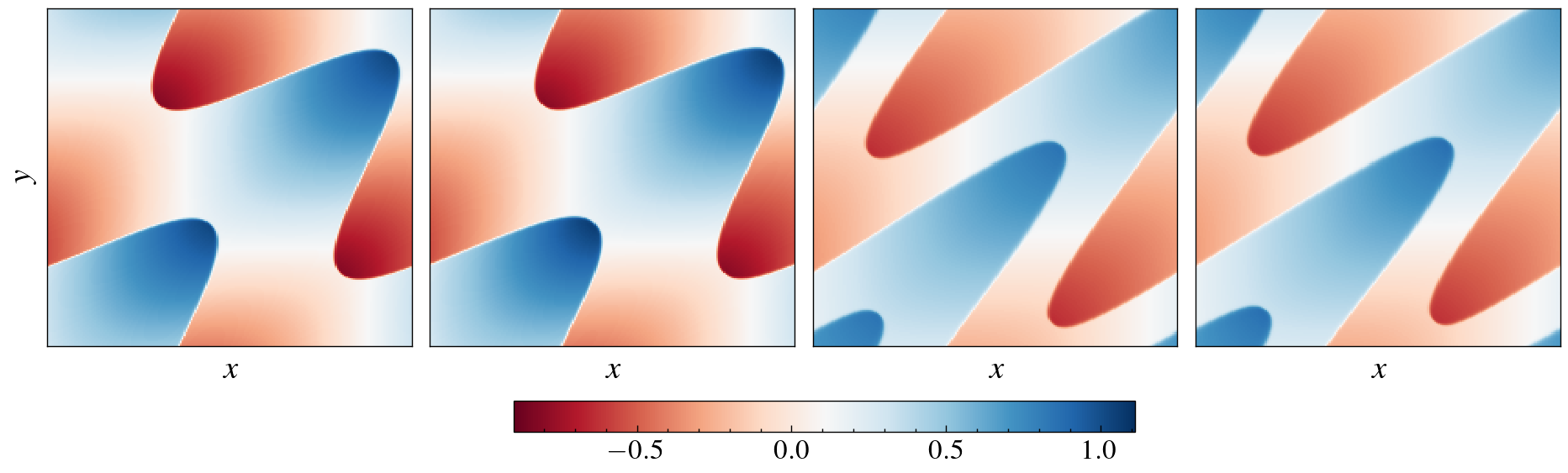}
    \caption{Comparison of the reference solution  of  \(u\) in 2D Burgers' equation \cref{eq: Burgers' equation 2d} with KT-enhanced CFN predictions for noise coefficient $\eta = 1$: (left) reference solution at t = .4, (middle-left) predictions at t = .4, (middle-right) reference solution at t = .8, (right) predictions at t = .8.} 
    \label{fig: 2d Burgers Noise}
\end{figure}
The noise introduced to the training data for the 2D Burgers' equation is analogous to that in \cref{eq: noise burgers}, and we again consider \(\eta = 0, .25, .5, 1\). Recall also that $L = L_{train} = 20$ and $N_{traj} = 5$.  \Cref{fig: 2d Burgers Noise} compares  the reference solutions with the KT-enhanced CFN predictions for $t = .4$ and $.8$. The predictions appear accurate even when $\eta = 1$.  
\begin{figure}[h!]
    \centering
    \includegraphics[width=0.9\linewidth]{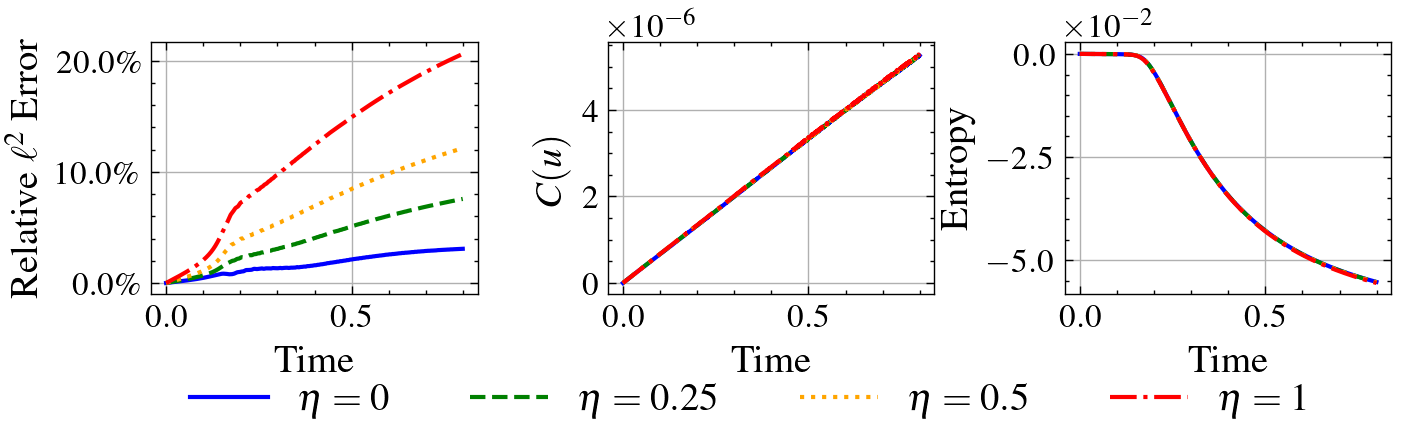}
    \caption{(left) Relative $\ell^2$ prediction error \cref{eq:relative l2 error}, (middle) Discrete conserved quantity remainder $C(u)$ in \cref{eq:conservemetric}, (right) discrete entropy remainder $\mathcal{J}(u)$ in \cref{eq: discrete entropy} for 2D Burgers' equation, $t \in [0,.8]$, and $\eta = 0, .25, .5, 1$ in \cref{eq: noise burgers}. }
    \label{fig: 2d Burgers Noise Error}
\end{figure}
 \Cref{fig: 2d Burgers Noise Error} displays the  relative \(\ell^2\) prediction error for $t \in [0,.8]$, once again confirming that the prediction maintains consistency for various noise levels $\eta \in [0,1]$. The figure also shows the discrete conserved quantity remainder metric, \(C(u)\) in \cref{eq:conservemetric}, and the discrete entropy remainder, $\mathcal{J}(u)$ in \cref{eq: discrete entropy}. Even with only $5$ training trajectories,  the KT-enhanced CFN successfully preserves conservation and satisfies  non-positivity for the entropy condition. 

\subsection{Coarse Observations}
\label{subsec: coarse}
We now assume the observations are noise-free, i.e.~$\eta = 0$, but may be coarse.  In particular, the prototype conservation laws are trained at varying spatial grid sizes, with $n_{train} = 64, 128, 256,512$ for each 1D problem and $n_{train} = 25\times25, 50\times50, 100\times100, 200\times200$ for the 2D Burgers' equation experiment. The temporal step size $\Delta t$ is fixed so as not to influence the overall error and to ensure that the CFL condition is satisfied. We also fix the testing grid size with  \(n_{test} = 512\) for each 1D problem and  \(n_{test} = 200\times200\) for the 2D Burgers' equation experiment. Observe that \(\Delta x\) for the testing and training data must be accordingly adjusted for the KT-enhanced CFN algorithm.

\subsubsection{Burgers' Equation} 
\label{subsubsec: 1D Burgers Grids}

\begin{figure}[h!]
    \centering
    \includegraphics[width=0.9\linewidth]{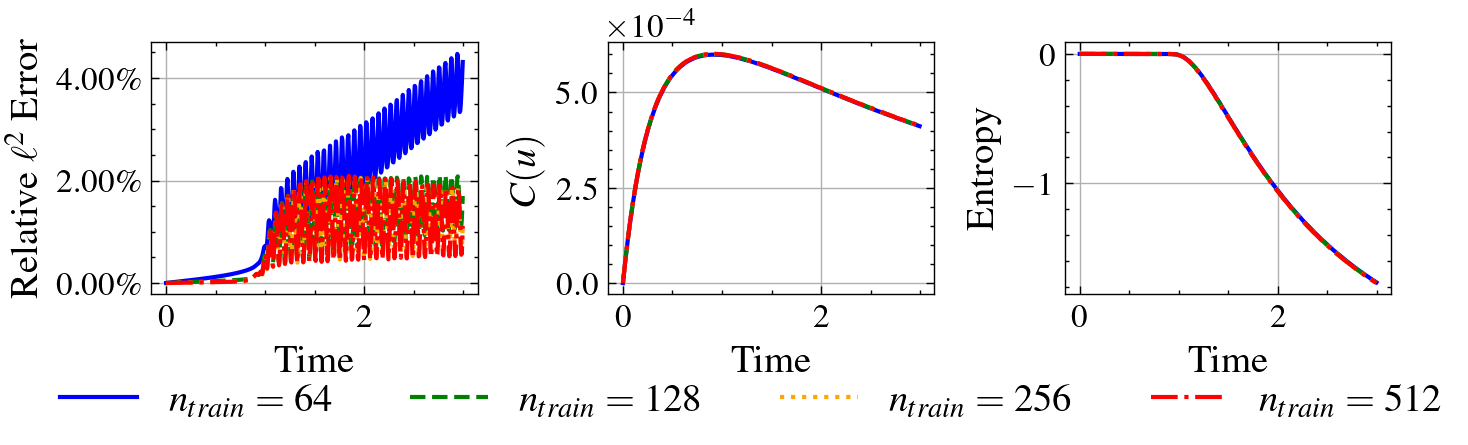}
    \caption{(left) Relative $\ell^2$ prediction error \cref{eq:relative l2 error}, (middle) Discrete conserved quantity remainder $C(u)$ in \cref{eq:conservemetric}, (right) discrete entropy remainder $\mathcal{J}(u)$ in \cref{eq: discrete entropy} for 1D Burgers' equation, $t \in [0,3]$, and $n_{train} = 64, 128, 256, 512$. }
    \label{fig: Burgers Grids Error}
\end{figure}

\Cref{fig: Burgers Grids Error} (left) compares the relative \(\ell^2\) prediction error for varying training spatial resolutions over the temporal domain $t \in [0,3]$, while \Cref{fig: Burgers Grids Error} (middle and right) displays the corresponding evolution of the discrete conserved quantity remainder \(C(u)\) and the discrete entropy remainder $\mathcal{J}(u)$  of the  KT-enhanced CFN predictions. Note the similarities to the behaviors observed for various noise coefficients $\eta$ in  \Cref{subsubsec: 1D Burgers Noise}. Predictions at different times exhibit the same general behavioral patterns as seen in \Cref{fig: Burgers Noise} and are therefore not shown.
\subsubsection{Shallow Water Equations}
\label{subsubsec: Shallow Water Grids}
\begin{figure}[h!]
    \centering
    \includegraphics[width=1.\linewidth]{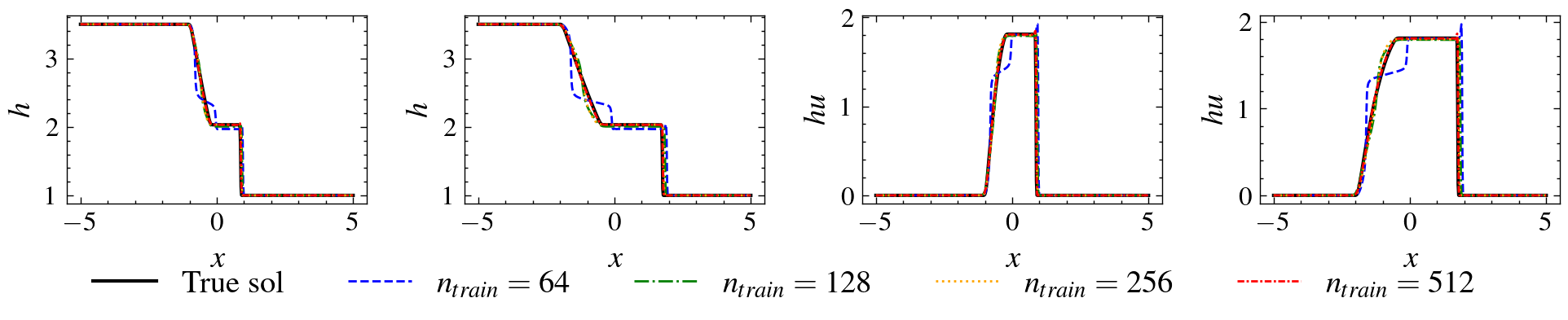}
     \caption{Comparison of the reference solution (black solid line) of height \(h\) and momentum \(hu\) for the shallow water equations with the  KT-enhanced CFN predictions for $n_{train} = 64, 128, 256, 512$: (left) \(t = .5\) of \(h\), (middle-left) \( t= 1\) of \(h\), (middle-right) \(t = .5\) of \(hu\), (right) \(t = 1\) of \(hu\).}
     \label{fig: Shallow Water Grids h}
\end{figure}

\begin{figure}[h!]
    \centering
    \includegraphics[width=1.\linewidth]{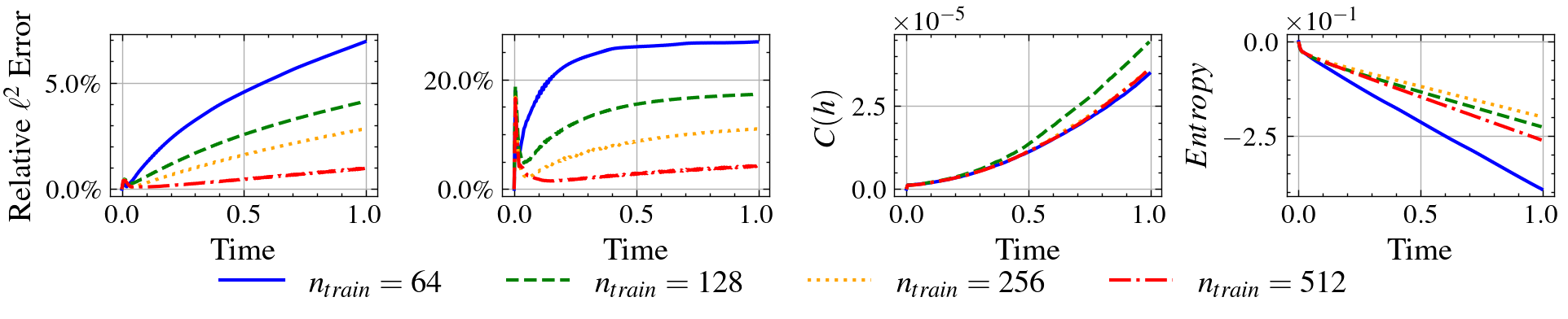}
 \caption{Relative \(\ell^2\) prediction error \cref{eq:relative l2 error} in the shallow water equations: (left) \(h\) (middle-left) \(hu\), and (middle-right) discrete conserved quantity remainder \cref{eq:conservemetric} of height $C(h)$; (right) discrete entropy remainder \(\mathcal{J}([h,hu]^T)\) in \cref{eq: discrete entropy} for the shallow water equations with $n_{train} = 64,128,256,512$.}
    \label{fig: Shallow Water Grids Error}
\end{figure}

\Cref{fig: Shallow Water Grids h} compares the predicted solution for height $h$ and momentum $hu$ at different times for varying $n_{train}$ values.  It is evident that having under-resolved training data can significantly degrade the accuracy of the solution, potentially introducing non-physical solutions with incorrect shock locations. Nevertheless, \Cref{fig: Shallow Water Grids Error} demonstrates convergence of the relative $\ell^2$ prediction errors for \(h\) and \(hu\) with respect to the training data resolution throughout the testing domain $[0,3]$. Although more analysis is needed to fully understand the training resolution requirements, it is promising that the KT-enhanced CFN  appear to follow  conventional numerical conservation law convergence properties with respect to resolution. Moreover, it may be possible to use the behavior observed in \Cref{fig: Shallow Water Grids Error} (left and middle-left) to determine if the resolution requirements are met. 


Finally \Cref{fig: Shallow Water Grids Error} (middle-right and right) displays the discrete conserved quantity remainders \(C({h})\) from \cref{eq:conservemetric} along with the discrete entropy remainder $\mathcal{J}([h, hu]^T)$ in \cref{eq: discrete entropy} for the shallow water equations for different choices of $n_{train}$.  We observe that conservation is preserved within $\mathcal{O}(10^{-5})$ and that the discrete entropy remainder is non-positive for all choices of $n_{train}$.  This is in spite of the solution for under-resolved training data ($n_{train} = 64$) not properly capturing all structural information. 

\subsubsection{Euler's Equation}
\label{subsubsec: Euler Grids}
\begin{figure}[h!]
    \centering
    \includegraphics[width=1.\linewidth]{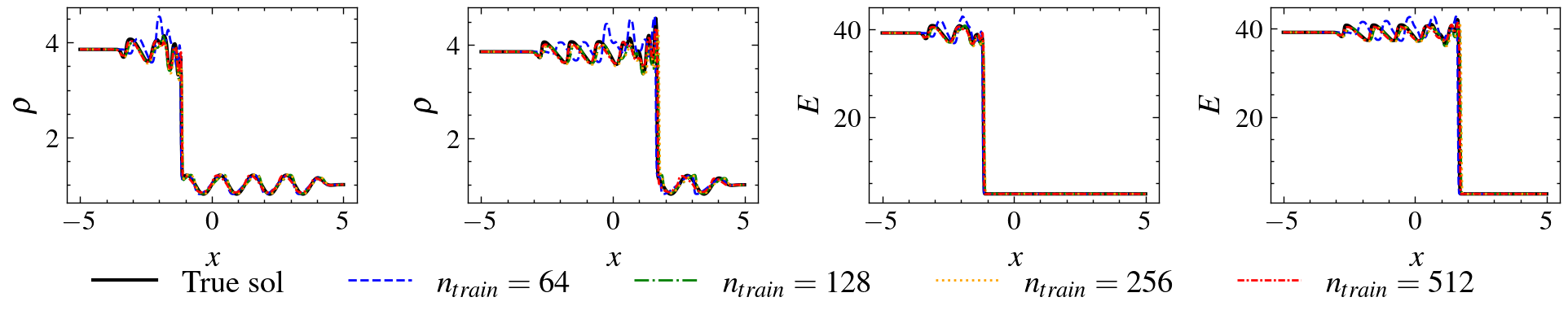}
    \caption{Comparison of the reference solution (black solid line) of density \(\rho\) and ennergy \(E\) in Euler's equation with the  KT-enhanced CFN predictions, $n_{train} = 64,128,256, 512$: (left) \(t = .8\) of \(\rho\) (middle-left) \(t = 1.6\) of \(\rho\) (middle-right) \(t = .8\) of \(E\) (right) \(t = 1.6\) of \(E\). }
    \label{fig: Euler Grids rho}
\end{figure}
\begin{figure}[h!]
    \centering
    \includegraphics[width=1.\linewidth]{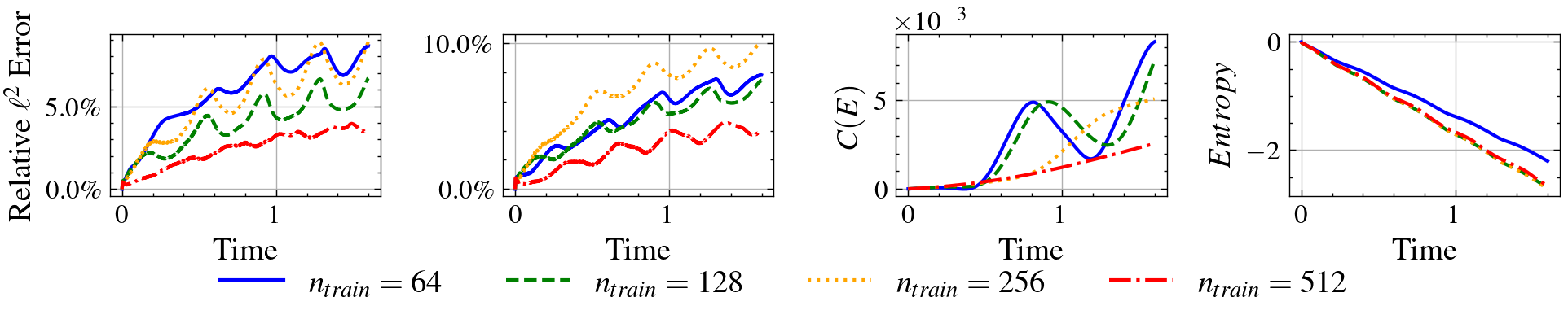}
    \caption{Relative $\ell_2$ prediction error \cref{eq:relative l2 error}: (left) \(\rho\), (middle-left) \(E\), and (middle-right) discrete conserved quantity remainder \(C({E})\) of \cref{eq:conservemetric}, and (right) discrete entropy remainder $\mathcal{J}([\rho, \rho u, E]^T)$ in \cref{eq: discrete entropy} for Euler's equation,  $n_{train} = 64,128, 256, 512$.}
    \label{fig: Euler Grids Error}
\end{figure}

\Cref{fig: Euler Grids rho} displays the KT-enhanced CFN predictions at different times for varying choices of $n_{train}$.\footnote{As in \Cref{subsubsec: Euler Noise} we omit the corresponding figure for momentum $\rho u$ as the results are comparable to those for energy $E$.} As was the case for the shallow water equations (see  \Cref{fig: Shallow Water Grids h}), it is evident that Euler's equations also have training resolution requirements for the KT-enhanced CFN.  While the prediction made with \(n_{train}=64\) clearly deviates the most from the solution for each time $t = .8$ and $1.6$, the relative \(\ell^2\) prediction errors \cref{eq:relative l2 error} shown in \Cref{fig: Euler Grids Error} (left and middle-left) suggest it yields an overall better performance. The relative \(\ell^2\) prediction error for momentum \(\rho u\) is omitted as the results are  similar to those for energy \(E\).  In this regard we first note that the scale of the error is smaller than that in the shallow water equation case.  Moreover, although the prediction using $n_{train} = 64$ does not capture the oscillatory behavior in smooth regions, the discretized {\em pointwise} error may be small.  By contrast in the $n_{train} = 256$ case, the largest error seems to occur as a result of the mismatch in shock {\em location}, and this apparently has the most significant impact on the overall error.  Thus we see that the relative \(\ell^2\) error, which is {\em global},  may not capture the convergence properties of the enhanced-CFN method (see \Cref{rem:relative l2 error}). \Cref{fig: Euler Grids Error} (middle-right and right) displays the discrete conserved quantity remainders \(C({E})\) with similar results for \(C({\rho})\) and \(C({\rho u})\), along with the discrete entropy remainder \(\mathcal{J}([\rho, \rho u, E]^T)\) i.e., \cref{eq: discrete entropy} for Euler's equation. 

\subsubsection{2D Burgers' Equation} \label{subsec:2DBurgunder-resolved}
\begin{figure}[h!]
    \centering
    \includegraphics[width=0.9\linewidth]{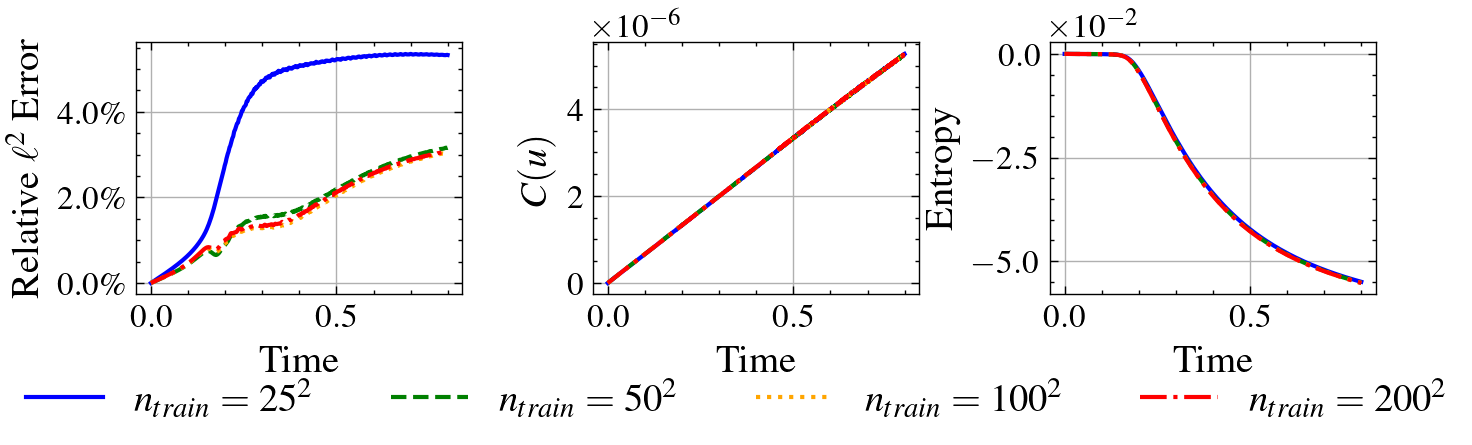}
    \caption{(left) Relative $\ell^2$ prediction error \cref{eq:relative l2 error}, (middle) Discrete conserved quantity remainder $C(u)$ in \cref{eq:conservemetric}, (right) discrete entropy remainder $\mathcal{J}(u)$ in \cref{eq: discrete entropy} for 2D Burgers' equation, $t \in [0,.8]$, and $n_{train} = 25^2, 50^2, 100^2, 200^2$. }
    \label{fig: 2d Burgers Grids Error}
\end{figure}

\Cref{fig: 2d Burgers Grids Error} demonstrates the impact of the choice of $n_{train}$ on the 2D Burgers' equation.  Similar to the 1D case, we observe that significant error results from insufficient resolution of the training data.  Once the resolution requirement is met (here $n_{train} = 50\times 50$) there does not seem to be any increased benefits in training on more data.  Moreover, conservation is preserved even in the under-resolved training environment, and the entropy condition is similarly met.   The KT-enhanced CFN predictions for different training resolutions follow the same behavior pattern as what is seen in \Cref{fig: 2d Burgers Noise} and is therefore omitted here.

\subsection{Sensitivity Analysis}
\label{subsec: sensitivity}

We now analyze the sensitivity of the entropy-stable CFN approach in the context of its specifications. In particular, we analyze the impacts of (1) time integration, (2) the effect of using a slope limiter, and (3) the underlying scheme's spatial order of accuracy. Through a number of numerical experiments (not all reported here) we conclude that the {\em combination of} (1) TVDRK3 time discretization, (2) the minmod slope limiter and (3)  spatial second order of accuracy ensures the  satisfactory prediction results  for varying levels of noise in the training data (\Cref{subsec: noisy}) as well for varying resolution (\Cref{subsec: coarse}).   We now extend this analysis to gain more insight on the importance of these underlying numerical properties.  We will also discuss the impact of the spectral radius approximation in \cref{eq: rho_W}. The experiments described below are conducted with the same parameters as those in \Cref{subsec: noisy} with noise coefficient $\eta = 1$. 

\subsubsection*{Time integration for stability and accuracy}
\begin{figure}[h!]
    \centering
    \includegraphics[width=0.7\linewidth]{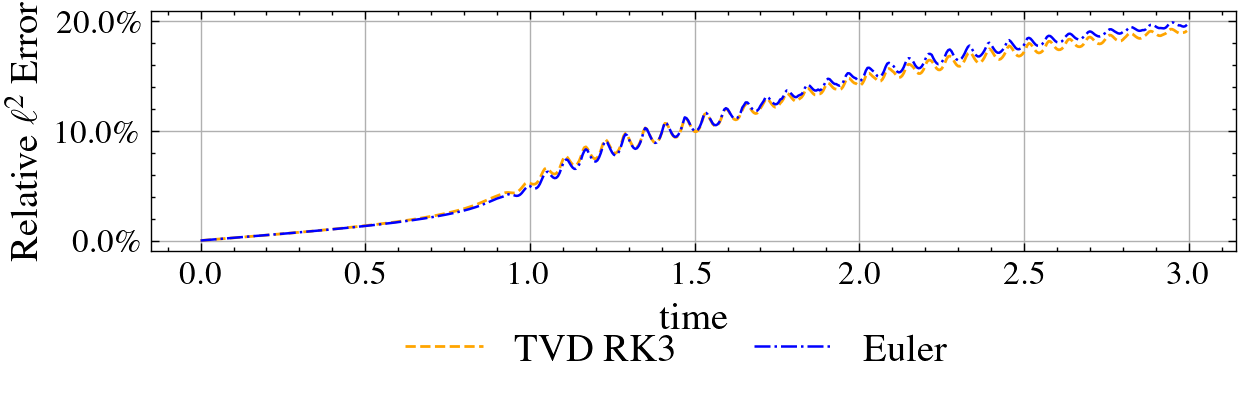}
    \caption{Relative \(\ell^2\) prediction error \cref{eq:relative l2 error} for $t \in [0,3]$ for Euler and TVDRK3  KT-enhanced CFN for 1D Burgers' equation with noisy training data \cref{eq: noise burgers}. Here $\eta = 1$.}
    \label{fig: 1d Burgers euler scheme}
\end{figure}

We chose TVDRK3 for temporal integration in our experiments because it is known to be stable for classical numerical conservation laws \cite{LeVeque92}. By contrast, it is also well known that forward Euler's method is not stable, which is easily verified through eigenvalue analysis \cite{LeVeque92}.  However, \Cref{fig: 1d Burgers euler scheme} suggests that  both time integration techniques  yield {\em stable}  KT-enhanced CFN predictions for 1D Burgers equation with noise coefficient \(\eta = 1\) in the training data \cref{eq: noise burgers}. Although beyond the scope of this investigation, we speculate that the recurrent loss function \cref{eq:recurrentloss} contributes enough dissipation to the predicted solution so that the forward Euler's method remains stable.  We note  that TVDRK3 yields about \(1\%\) better {\em accuracy} over time. 

\subsubsection*{Spatial order of accuracy and slope limiting}
\begin{figure}[h!]
    \centering
    \includegraphics[width=0.9\linewidth]{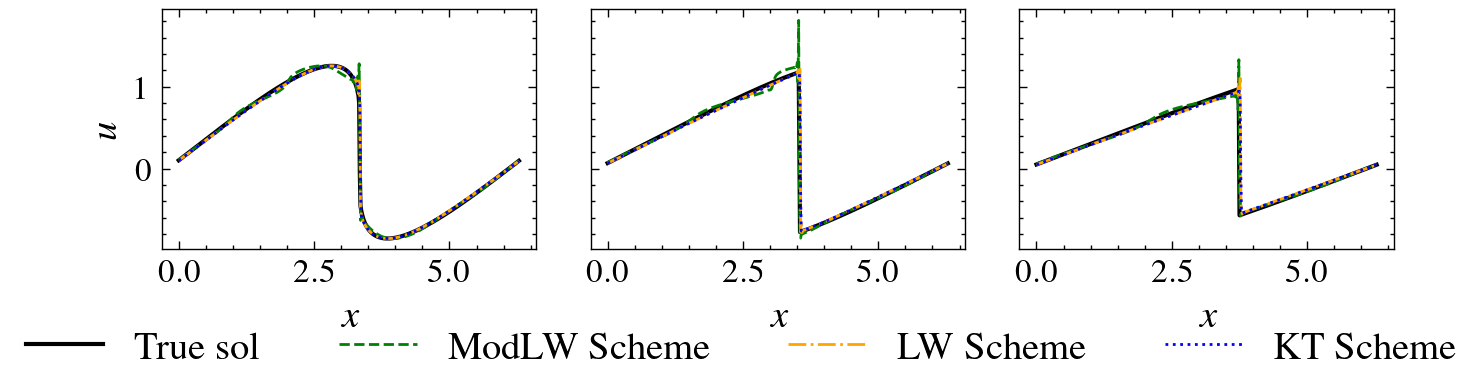}
    \caption{Comparison of the reference solution (black solid line) to 1D Burgers' equation with predictions from the KT-enhanced CFN \cref{eq: Numerical flux NN}, the LW-enhanced CFN \cref{eq: LWS-cfn}, and the modLW-enhanced CFN \cref{eq: modLWS-cfn} at (left) \(t = 1\) (middle) \(t = 2\) (right) \(t = 3\) for $\eta = 1$ in \cref{eq: noise burgers}.}
    \label{fig: Burgers Scheme}
\end{figure}

\Cref{fig: Burgers Scheme} 
compares the performance of the KT-enhanced CFN \cref{eq: Numerical flux NN}, the LW-enhanced CFN \cref{eq: LWS-cfn} and modLW-enhanced CFN \cref{eq: modLWS-cfn} for 1D Burgers equation where $\eta = 1$ in \cref{eq: noise burgers}.  In all cases we used TVDRK3 for temporal integration.

\Cref{fig: Burgers Scheme} demonstrates that both the LW-enhanced CFN and modLW-enhanced CFN predictions introduce artificial oscillations as the shock forms, which become more pronounced in time. Since all three methods are designed to have second order spatial accuracy in classical numerical conservation laws, we can conclude that, as expected, spatial order of accuracy does not prevent the introduction of artificial oscillations in the CFN solution.  In this regard it appears that the slope limiter used in the KT-enhanced CFN is crucial for mitigating the  spurious oscillations and overshoots in the prediction.  Finally we note that the shock characteristic structure is maintained for long term prediction for each method, which is a feature of higher order accuracy numerical methods.\footnote{This is assuming, of course, that the method is written in proper conservative form.} 

\begin{figure}[h!]
    \centering
    \includegraphics[width=0.9\linewidth]{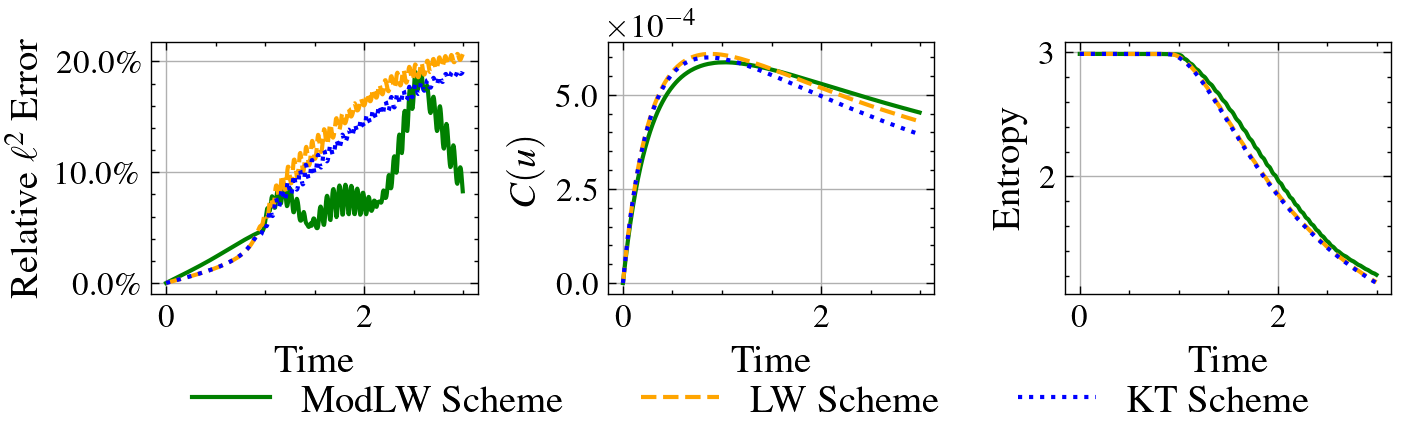}
    \caption{(left) Relative \(\ell^2\) prediction error \cref{eq:relative l2 error}, and (middle) discrete conserved quantity remainder \(C(u)\) in \cref{eq:conservemetric} and  (right) discrete entropy remainder  \(\mathcal{J}(u)\) in \cref{eq: discrete entropy} for 1D Burgers' equation given noisy training data \cref{eq: noise burgers} with \(\eta=1\) in $t = [0,3]$ for the modLW-enhanced CFN, the LW-enhanced CFN, and the KT-enhanced CFN predictions.}
    \label{fig: Burgers Scheme Error}
\end{figure}

There are other measurable impacts on the overall prediction performance based on underlying flux choice, however. \Cref{fig: Burgers Scheme Error} shows the relative \(\ell^2\) prediction error \cref{eq:relative l2 error} for the KT-enhanced CFN, the LW-enhanced CFN, and the modLW-enhanced CFN predictions in the time domain $[0,3]$. The accurate shock speed obtained using the modLW-enhanced CFN results in the smallest $\ell^2$ prediction error after the shock forms at $t = 1$.  This is in spite of  oscillatory behavior and a noticeable overshoot in the solution at the shock location, which also suggests caution when using the relative \(\ell^2\) prediction error as a metric in cases with discontinuities (see \Cref{rem:relative l2 error}). 
All three enhanced-CFN approaches ensure both conservation and entropy non-positivity are enforced, as shown in \Cref{fig: Burgers Scheme Error}. By contrast, the original CFN method \cite{chencfn} does not enforce  entropy non-positivity.  It also does not preserve the shock structure in the $\eta = 1$ case, even when trained over different training periods that allow the network to train on data that includes both smooth and discontinuous solution profiles.

\begin{rem}
\label{rem:initialconditions} 
It is important to emphasize that the enhanced CFN method is able to predict  future solutions from a limited temporal domain of training data that comes from an unknown underlying physical system, that is, without information regarding how the solution may evolve.  Essentially this means that the learning is extrapolatory. This is especially apparent in the Burgers equation experiments since the training data contain no profiles with discontinuities.  This has potential implications regarding the compressibility of training data for machine learning algorithms of physical problems admitting conservation properties.
\end{rem}

\subsubsection*{Computation of KT-enhanced CFN spectral radius} 
\begin{figure}[h!]
    \centering
    \includegraphics[width=1.\linewidth]{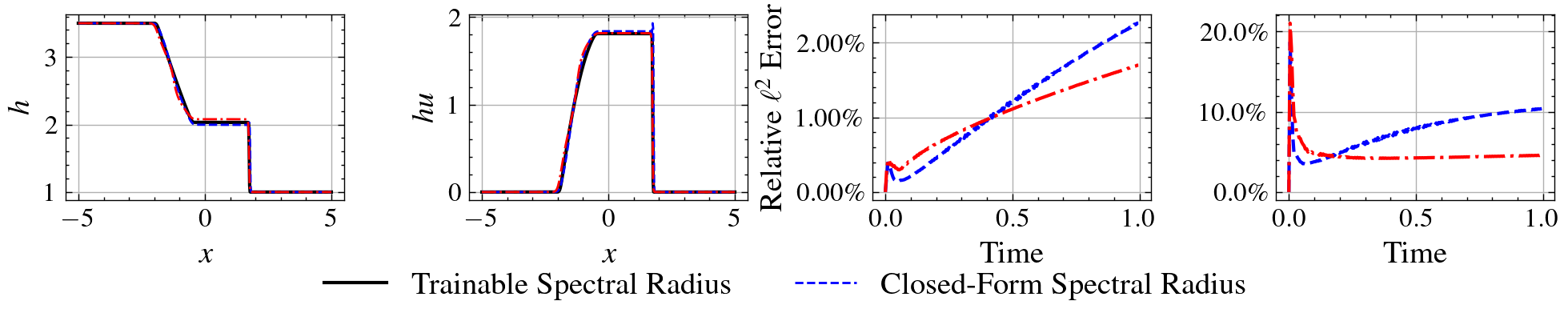}
    \caption{KT-enhanced CFN predictions at time $t = 1$ with closed form and trained spectral radius computation for the shallow water equations with $\eta = 1$ in \cref{eq:noise_shallow} (left) $h$, (middle-left) $hu$. Relative $\ell^2$ prediction error \cref{eq:relative l2 error} in testing domain $[0,1]$. (middle-right) height $h$.  (right) momentum $hu$.}
    \label{fig: shallow water closed form}
\end{figure}

Recall that \cref{eq: Numerical flux NN} requires an approximation to the maximum wave speed \cref{eq: Maximum wave speed NN}.  As noted there, determining the spectral radius of the corresponding Jacobian matrix is inherently problem dependent.  Since the flux is unknown,  we approximate the maximum wave speed using \cref{eq: rho_W}, which relies on a {\em learned}  spectral radius.  \Cref{fig: shallow water closed form} shows the impact of this approximation, specifically by  comparing the KT-enhanced CFN predictions for the shallow water equations to those that would have a closed-form spectral radius available for computing  \cref{eq: Maximum wave speed NN}.   The relative \(\ell^2\) prediction error for both cases is also shown.  Although using \cref{eq: rho_W} introduces some error, it does not appear to play a significant role in the overall accuracy.       

%% file: 6_summary.tex
\section{Concluding remarks}
\label{sec:conclusion}

This paper introduced the {\em entropy-stable} conservative flux form neural network (CFN) to predict solutions of conservation laws that admit discontinuous profiles.   Specifically, we incorporated the entropy-stable second-order non-oscillatory Kurganov-Tadmor (KT) scheme into the design of the neural networks. The method also uses a neural network to approximate the spectral   radius of Jacobian matrix in the KT scheme, which is simultaneously trained with the flux network.

Our numerical experiments  highlight the benefits of the entropy-stable CFN framework in capturing complex dynamics while maintaining stability and accuracy in both noisy and sparse observation environments.  This is even true when  the training domain does not itself  contain discontinuous solution profiles.  That is, the entropy-stable CFN is extrapolatory, and does not rely on oracle knowledge of possible long time solution profiles. Our experiments also demonstrate conservation and entropy stability over the testing period. The entropy-stable CFN method is an extension of the CFN method introduced in \cite{chencfn}, which used a first order finite volume conservative flux form and did not consider entropy stability in its design.  As such, it required a longer training period and could not effectively mitigate noise in all cases to obtain an accurate prediction.   

 Our investigation also  analyzed the influence of various specification options, for example in comparing the  KT-enhanced CFN with the CFN enhanced by other entropy stable  methods such as the Lax-Wendroff and modified Lax-Wendroff schemes. We observe that the non-oscillatory  slope-limiter in combination with the second order accuracy are critical in obtaining stable long term predictions.

There are several future directions for this work. One natural extension is to use techniques such as the discontinuous Galerkin (DG) method for higher-dimensional systems, especially given more complex geometries. Another avenue lies in refining the approximation of the local Lax–Friedrichs flux. A more precise method for approximating the monotone numerical flux may alleviate the need for computing the spectral radius with a neural network. Furthermore, the current study focuses on solving idealized hyperbolic conservation laws, each of which assumes access to full state variables. In practice, only observables like mass, density, velocity, and pressure may be measurable, however. This raises the important research question of how to {\em infer} underlying conservation laws from observable data. Additionally, since the entropy inequality shares a structural similarity with conservation laws, developing a parallel framework to learn entropy pairs should also be possible. Lastly, while our findings demonstrate the empirical success of the framework, a rigorous theoretical analysis of its convergence and stability properties remains an important area for future investigations.